\documentclass[11pt,a4paper]{article}
\usepackage[top=1in, bottom=1.25in, left=1in, right=1in]{geometry}
\usepackage[utf8]{inputenc}
\usepackage[english]{babel}
\usepackage{graphicx,epstopdf}
\usepackage{caption, subcaption,float}
\usepackage{color}
\usepackage{multirow}
\usepackage{tikz}
\usepackage{amsthm}
\usepackage{tabls}
\usepackage{booktabs}
\usepackage{cite}
\usepackage[round,comma,sort,longnamesfirst]{natbib}

\usepackage[ruled, linesnumberedhidden,vlined]{algorithm2e} 
\usetikzlibrary{matrix}
\usepackage{amsmath,amssymb,dsfont,mathtools}
\usepackage{hyperref}
\hypersetup{
    colorlinks=true,       
    linkcolor=blue,          
    citecolor=olive,       
    filecolor=magenta,      
    urlcolor=blue          
}
\usepackage{aliascnt}
\usepackage{titlesec}
\usepackage{soul}
\usepackage{transparent}

\begin{document}

\title{Fundamental techniques in the study of \\ parabolic subgroups of Artin groups}
\date{\today }
\author{Mar\'{i}a Cumplido}


\maketitle
\theoremstyle{plain}
\newtheorem{theorem}{Theorem}

\theoremstyle{plain}
\newtheorem{question}{Question}

\newaliascnt{lemma}{theorem}
\newtheorem{lemma}[lemma]{Lemma}
\aliascntresetthe{lemma}
\providecommand*{\lemmaautorefname}{Lemma}

\newaliascnt{proposition}{theorem}
\newtheorem{proposition}[proposition]{Proposition}
\aliascntresetthe{proposition}
\providecommand*{\propositionautorefname}{Proposition}

\newaliascnt{corollary}{theorem}
\newtheorem{corollary}[corollary]{Corollary}
\aliascntresetthe{corollary}
\providecommand*{\corollaryautorefname}{Corollary}

\newaliascnt{conjecture}{theorem}
\newtheorem{conjecture}[conjecture]{Conjecture}
\aliascntresetthe{conjecture}
\providecommand*{\conjectureautorefname}{Conjecture}

\newtheorem{exercise}{Exercise}
\providecommand*{\exerciseautorefname}{Exercise}

\theoremstyle{remark}

\newaliascnt{claim}{theorem}
\newaliascnt{remark}{theorem}

\newtheorem{claim}[claim]{Claim}
\newtheorem{remark}[remark]{Remark}
\newaliascnt{notation}{theorem}
\newtheorem{notation}[notation]{Notation}
\aliascntresetthe{notation}
\providecommand*{\notationautorefname}{Notation}

\aliascntresetthe{claim}
\providecommand*{\claimautorefname}{Claim}

\aliascntresetthe{remark}
\providecommand*{\remarkautorefname}{Remark}

\newtheorem*{claim*}{Claim}
\theoremstyle{definition}

\newaliascnt{definition}{theorem}
\newtheorem{definition}[definition]{Definition}
\aliascntresetthe{definition}
\providecommand*{\definitionautorefname}{Definition}

\newaliascnt{example}{theorem}
\newtheorem{example}[example]{Example}
\aliascntresetthe{example}
\providecommand*{\exampleautorefname}{Example}

\newenvironment{solution}
  {\begin{proof}[Solution]}
  {\end{proof}}
  
\titleformat{\paragraph}[block]{\normalfont\normalsize\bfseries}{\theparagraph}{1em}{}

\def\autorefspace{\hspace*{-0.5pt}}
\def\sectionautorefname{Section\autorefspace}
\def\subsectionautorefname{Section\autorefspace}
\def\subsubsectionautorefname{Section\autorefspace}
\def\figureautorefname{Figure\autorefspace}
\def\subfigureautorefname{Figure\autorefspace}
\def\tableautorefname{Table\autorefspace}
\def\equationautorefname{Equation\autorefspace}
\def\Itemautorefname{item\autorefspace}
\def\Hfootnoteautorefname{footnote\autorefspace}
\def\AMSautorefname{Equation\autorefspace}

\newcommand{\co}{\simeq_c}
\newcommand{\w}{\widetilde}
\newcommand{\po}{\preccurlyeq}
\newcommand{\soo}{\succcurlyeq}
\newcommand{\dist}{\mathrm{d}}

\def\Z{\mathbb Z} 
\def\Ker{{\rm Ker}} \def\R{\mathbb R} \def\GL{{\rm GL}}
\def\HH{\mathcal H} \def\C{\mathbb C} \def\P{\mathbb P}
\def\SSS{\mathfrak S} \def\BB{\mathcal B} \def\PP{\mathcal P} 
\def\supp{{\rm supp}} \def\Id{{\rm Id}} \def\Im{{\rm Im}}
\def\MM{\mathcal M} \def\S{\mathbb S}
\newcommand{\bigveer}{\bigvee^\Lsh}
\newcommand{\wedger}{\wedge^\Lsh}
\newcommand{\veer}{\vee^\Lsh}
\def\diam{{\rm diam}}

\newcommand{\myref}[2]{\hyperref[#1]{#2~\ref*{#1}}}

\begin{abstract}
This survey was written on the occasion of the course I gave at the Winterbraids XIV workshop in Bordeaux (2025). Its main purpose is to present the techniques that have proven most effective in the study of parabolic subgroups of Artin groups, with particular emphasis on the parabolic subgroups intersection problem. The survey highlights the core ideas and strategies behind them, aiming to give the reader a concise and accessible entry point to the essential methods.

\medskip

{\footnotesize

}

\end{abstract}

\noindent\textbf{Note to the reader.}
This survey brings together several results and ideas due to other researchers, most notably Yago Antolín, Martín Blufstein, Ruth Charney, Islam Foniqi, Federica Gavazzi, Volker Gebhardt, Eddy Godelle, Juan González-Meneses, Ivan Marin, Alexandre Martin, Philip Möller, Rose Morris-Wright, Luis Paris, Olga Varghese, Nikolas Vaskou, and Bert Wiest. Along the way, I have added explanations that would not usually appear in a research paper, with the aim of weaving these contributions into a broader picture of the recent history of parabolic subgroups. To help the reader engage with the material, I have also included a series of introductory-level exercises. The solutions are short and straightforward, so I encourage beginning readers to attempt them. Complete solutions can be found at the end of the survey.

\section{Introduction}

Braid groups, introduced by \cite{Artin1947}, lie at the intersection of algebra and topology. Their connections with low-dimensional topology—such as mapping class groups and knot theory—as well as with various algebraic and combinatorial structures—like Garside theory, algorithmics, cryptography, and category theory—have attracted sustained interest over the past 75 years.

On the one hand, the braid group on $n$ strands, denoted $\mathcal{B}_n$, can be defined as the mapping class group of the $n$-punctured disc $\mathcal{D}_n$, that is, the group of isotopy classes of orientation-preserving homeomorphisms of $\mathcal{D}_n$ that fix the boundary pointwise. On the other hand, it admits the following algebraic presentation:
\[
\mathcal{B}_n = \left\langle \sigma_1, \dots, \sigma_{n-1} \,\middle|\, 
\begin{array}{ll}
\sigma_i \sigma_j = \sigma_j \sigma_i & \text{if } |i - j| > 1, \\
\sigma_i \sigma_j \sigma_i = \sigma_j \sigma_i \sigma_j & \text{if } |i - j| = 1
\end{array}
\right\rangle.
\]

\begin{figure}[H]
  \centering
  \def\svgwidth{0.60\columnwidth} 
  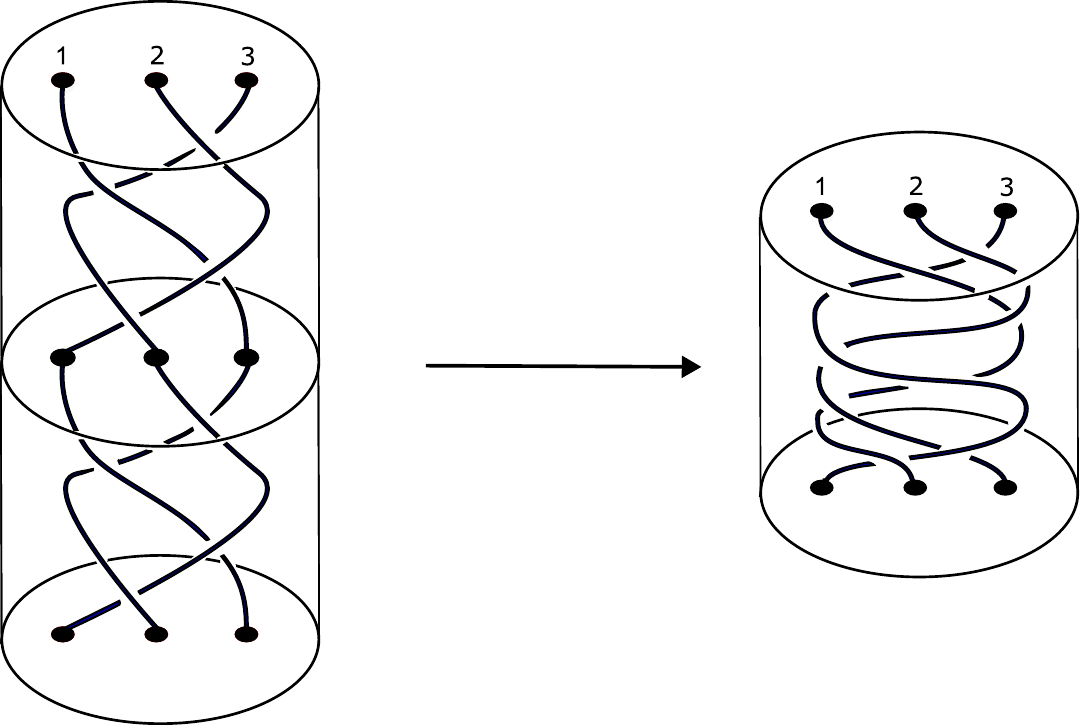
  \caption{Multiplication of two braids.}\label{Dibujo:multiplicar}
\end{figure}

As we will explain, these two perspectives coincide with the geometric notion of a \emph{braid} on $n$ strands: an object embedded in a cylinder, with $n$ distinguished points on both the top and bottom discs. Each top point is connected to a unique bottom point by a path (called a \emph{strand}) that monotonically runs in the vertical direction, and such that all strands are pairwise disjoint. Two braids are considered equivalent if one can be continuously deformed into the other through such configurations, without allowing strands to intersect. The group operation is defined by stacking one braid on top of another and rescaling the vertical direction, as illustrated in \autoref{Dibujo:multiplicar}.

\begin{figure}[H]
  \centering
  \def\svgwidth{0.7\columnwidth} 
\begingroup%
  \makeatletter%
  \providecommand\color[2][]{%
    \errmessage{(Inkscape) Color is used for the text in Inkscape, but the package 'color.sty' is not loaded}%
    \renewcommand\color[2][]{}%
  }%
  \providecommand\transparent[1]{%
    \errmessage{(Inkscape) Transparency is used (non-zero) for the text in Inkscape, but the package 'transparent.sty' is not loaded}%
    \renewcommand\transparent[1]{}%
  }%
  \providecommand\rotatebox[2]{#2}%
  \newcommand*\fsize{\dimexpr\f@size pt\relax}%
  \newcommand*\lineheight[1]{\fontsize{\fsize}{#1\fsize}\selectfont}%
  \ifx\svgwidth\undefined%
    \setlength{\unitlength}{221.34700265bp}%
    \ifx\svgscale\undefined%
      \relax%
    \else%
      \setlength{\unitlength}{\unitlength * \real{\svgscale}}%
    \fi%
  \else%
    \setlength{\unitlength}{\svgwidth}%
  \fi%
  \global\let\svgwidth\undefined%
  \global\let\svgscale\undefined%
  \makeatother%
  \begin{picture}(1,0.40765196)%
    \lineheight{1}%
    \setlength\tabcolsep{0pt}%
    \put(0,0){\includegraphics[width=\unitlength,page=1]{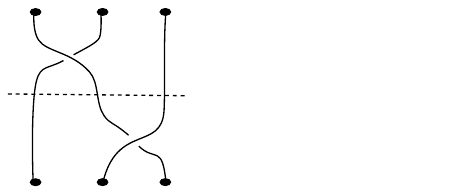}}%
    \put(0.00220635,0.27988103){\color[rgb]{0,0,0}\makebox(0,0)[lt]{\lineheight{1.25}\smash{\begin{tabular}[t]{l}$\sigma_1$\end{tabular}}}}%
    \put(0.71885956,0.20441788){\color[rgb]{0,0,0}\makebox(0,0)[lt]{\lineheight{1.25}\smash{\begin{tabular}[t]{l}$\sigma_1$\end{tabular}}}}%
    \put(0.00028565,0.10200616){\color[rgb]{0,0,0}\makebox(0,0)[lt]{\lineheight{1.25}\smash{\begin{tabular}[t]{l}$\sigma_2^{-1}$\end{tabular}}}}%
    \put(0,0){\includegraphics[width=\unitlength,page=2]{generadoresmcg.pdf}}%
    \put(0.84736716,0.20368377){\color[rgb]{0,0,0}\makebox(0,0)[lt]{\lineheight{1.25}\smash{\begin{tabular}[t]{l}$\sigma_2^{-1}$\end{tabular}}}}%
  \end{picture}%
\endgroup%

  \caption{How a braid decomposes into standard generators, and what these generators look like as mapping classes in the braid group on three strands.}\label{Dibujo:generadores}
\end{figure}

If we interpret the height of the cylinder as a time variable, a braid describes the motion of the $n$ points in the disc over time—precisely, an element of the mapping class group of $\mathcal{D}_n$. On the other hand, the \emph{standard generator} $\sigma_i$ corresponds to the braid obtained from the trivial braid—where all strands are vertical and no crossings occur—by introducing a single crossing between the strand in position $i$ and the one in position $i+1$. The inverse $\sigma_i^{-1}$ represents the same crossing, but in the opposite direction. It is easy to see that any braid can be expressed as a product of these generators and their inverses. Moreover, from the viewpoint of the mapping class group, each generator $\sigma_i$ corresponds to a half-twist exchanging the $i$-th and $(i+1)$-th punctures in~$\mathcal{D}_n$ (see \autoref{Dibujo:generadores}).

The study of braid groups has greatly benefited from both topological and combinatorial perspectives, providing rich tools and insights into their structure. A natural question, then, is whether these methods and properties extend to broader classes of groups. One such class is that of \emph{Artin--Tits groups} (or simply \emph{Artin groups}), introduced by Jacques Tits in the 1960s. These groups arise in geometric contexts, such as hyperplane arrangements and Coxeter theory, yet they are defined through a simple combinatorial presentation that generalizes that of braid groups.

To define an Artin group, we start with a finite set $S$, and for each pair $s, t \in S$ we assign a label $m_{s,t} \in \{2, 3, \dots, \infty\}$. These labels encode all the relations in the presentation. Specifically, if $m_{s,t} \ne \infty$, we impose the relation:
\[
\underbrace{sts\cdots}_{m_{s,t} \text{ terms}} = \underbrace{tst\cdots}_{m_{s,t} \text{ terms}},
\]
where the words on both sides alternate between $s$ and $t$, and have length $m_{s,t}$. If $m_{s,t} = \infty$, then no relation is imposed between $s$ and $t$.
This data can also be encoded in a labelled complete graph, known as the \emph{Coxeter graph} $\Gamma$, whose vertices are the elements of $S$, and whose edges are labelled with the corresponding $m_{s,t}$. The group defined in this way is denoted:
\[
A[\Gamma] = \left\langle S \;\middle|\; \underbrace{sts\cdots}_{m_{s,t}} = \underbrace{tst\cdots}_{m_{s,t}} \;\text{ for all } s,t\in S \text{ with } m_{s,t} \ne \infty \right\rangle.
\]
When the context is clear, we may abuse notation and write $A_S$ instead of $A[\Gamma]$. When we do not state it explicitly, we will always assume that $S=\{s_1,s_2,s_3,\cdots \}$.

\medskip

Depending on the family of Artin groups under consideration, different conventions are commonly used to simplify the Coxeter graph:
\begin{itemize}
  \item In the \emph{no-$2$-convention}, edges  corresponding to $m_{s,t}=2$ are omitted and edges corresponding to $m_{s,t}=3$ are not labelled. This is useful when studying groups with many commutations, e.g. braid groups, and it is the convention used in Dynkin diagrams.
  \item In the \emph{no-$\infty$-convention}, edges labelled $\infty$ are omitted, which is useful in the study of groups with many $\infty$'s, e.g. right-angled Artin groups (RAAGs).
  \item In the \emph{all-edges convention}, the graph is complete and all labels $m_{s,t}$, including $2$ and $\infty$, are explicitly displayed.
\end{itemize}

In these notes, we will make use of all three conventions as appropriate. Each highlights different structural properties of the Artin group under study. We say that an Artin group is \emph{irreducible} if its Coxeter graph is connected (in the no-$2$-convention).

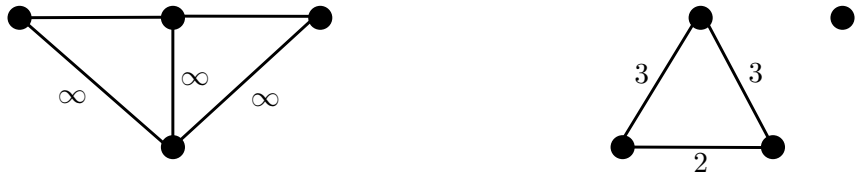
\begin{figure}[H]
  \centering
  \def\svgwidth{0.7\columnwidth} 
\begingroup%
  \makeatletter%
  \providecommand\color[2][]{%
    \errmessage{(Inkscape) Color is used for the text in Inkscape, but the package 'color.sty' is not loaded}%
    \renewcommand\color[2][]{}%
  }%
  \providecommand\transparent[1]{%
    \errmessage{(Inkscape) Transparency is used (non-zero) for the text in Inkscape, but the package 'transparent.sty' is not loaded}%
    \renewcommand\transparent[1]{}%
  }%
  \providecommand\rotatebox[2]{#2}%
  \newcommand*\fsize{\dimexpr\f@size pt\relax}%
  \newcommand*\lineheight[1]{\fontsize{\fsize}{#1\fsize}\selectfont}%
  \ifx\svgwidth\undefined%
    \setlength{\unitlength}{179.13514025bp}%
    \ifx\svgscale\undefined%
      \relax%
    \else%
      \setlength{\unitlength}{\unitlength * \real{\svgscale}}%
    \fi%
  \else%
    \setlength{\unitlength}{\svgwidth}%
  \fi%
  \global\let\svgwidth\undefined%
  \global\let\svgscale\undefined%
  \makeatother%
  \begin{picture}(1,0.1939397)%
    \lineheight{1}%
    \setlength\tabcolsep{0pt}%
    \put(0,0){\includegraphics[width=\unitlength,page=1]{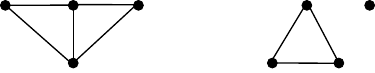}}%
    \put(0.74013822,0.10295829){\color[rgb]{0,0,0}\makebox(0,0)[lt]{\lineheight{1.25}\smash{\begin{tabular}[t]{l}$3$\end{tabular}}}}%
    \put(0.87507212,0.10327719){\color[rgb]{0,0,0}\makebox(0,0)[lt]{\lineheight{1.25}\smash{\begin{tabular}[t]{l}$3$\end{tabular}}}}%
    \put(0.81013059,-0.00321328){\color[rgb]{0,0,0}\makebox(0,0)[lt]{\lineheight{1.25}\smash{\begin{tabular}[t]{l}$2$\end{tabular}}}}%
    \put(0.28728362,0.07580438){\color[rgb]{0,0,0}\makebox(0,0)[lt]{\lineheight{1.25}\smash{\begin{tabular}[t]{l}$\infty$\end{tabular}}}}%
    \put(0.05924554,0.07824031){\color[rgb]{0,0,0}\makebox(0,0)[lt]{\lineheight{1.25}\smash{\begin{tabular}[t]{l}$\infty$\end{tabular}}}}%
    \put(0.20405192,0.09996843){\color[rgb]{0,0,0}\makebox(0,0)[lt]{\lineheight{1.25}\smash{\begin{tabular}[t]{l}$\infty$\end{tabular}}}}%
  \end{picture}%
\endgroup%

  \caption{The Coxeter graph corresponding to the free product of $\mathbb{Z}$ and a braid group on 3 strands. On the left with the no-$2$-convention and on the right with the no-$\infty$-convention.}\label{Dibujo:convenciones}
\end{figure}

From the same Coxeter graph~\(\Gamma\), one can also define the associated \emph{Coxeter group}~\(W[\Gamma]\), which has a similar presentation, but with all generators being involutions. Explicitly,
\[
W[\Gamma] = \left\langle S \;\middle|\; s^2 = 1 \text{ for all } s \in S, \quad \underbrace{sts\cdots}_{m_{s,t}} = \underbrace{tst\cdots}_{m_{s,t}} \;\text{ for all } s,t \in S \text{ with } m_{s,t} \ne \infty \right\rangle.
\]
When the set of relations is clear from context, we will also write~\(W_S\) instead of~\(W[\Gamma]\).

\medskip

There is a natural surjective homomorphism
\[
\theta \colon A[\Gamma] \longrightarrow W[\Gamma],
\]
which sends each Artin generator~\(s \in S\) to the corresponding Coxeter generator~\(s \in W[\Gamma]\). The kernel of $\theta$ is denoted by $\mathrm{CA}[\Gamma]$ and it is called the corresponding \emph{coloured} (or \emph{pure}) \emph{Artin group}. In the braid group on $n$ strands, the corresponding Coxeter group is the symmetric group on $n$ elements and the coloured Artin group is the pure braid group on $n$ strands, which consists of braids that induce the trivial permutation. A canonical set-theoretic section of~\(\theta\) is given by the inclusion
\[
\iota \colon W[\Gamma] \hookrightarrow A[\Gamma],
\]
which for each element~\(w \in W[\Gamma]\) chooses a word of minimal length and sends it to the corresponding element of $A[\Gamma]$ represented by this reduced word.

\medskip

We call a \emph{standard parabolic subgroup} (also known as a \emph{special subgroup}) any subgroup \( A_X \subseteq A_S \) of rank $k$ generated by a subset \( X \subseteq S \) such that $|X|=k$. By a classical result of ~\cite{Vanderlek}, the subgroup \( A_X \) is itself an Artin group associated with the full subgraph \( \Gamma_X \subseteq \Gamma \) spanned by the vertices in \( X \). In other words, \( A_X \cong A[\Gamma_X] \). More generally, any conjugate of a standard parabolic subgroup is called a \emph{parabolic subgroup}. That is, we say that a subgroup \( P \subseteq A_S \) is parabolic of rank $k$ if there exists \( \alpha \in A_S \) and \( X \subseteq S \), $|X|=k$, such that \( P = \alpha A_X \alpha^{-1} \). We say that a parabolic subgroup is \emph{irreducible} if the associated standard subgroup \( A_X \) is irreducible, i.e., if the subgraph \( \Gamma_X \) is connected under the no-2-convention. 

The same notions apply to Coxeter groups. Given a Coxeter group \( W_S \), a \emph{standard parabolic subgroup} is any subgroup \( W_X \subseteq W_S \) generated by a subset \( X \subseteq S \), and any conjugate \( w W_X w^{-1} \subseteq W_S \) is called a \emph{parabolic subgroup} of \( W_S \). 

Parabolic subgroups play a central role in the study of Artin groups, as they capture much of the internal structure of these groups. Understanding their behaviour—how they intersect, embed, and generate—is fundamental to developing a global understanding of Artin groups. The aim of this survey is to provide an overview of the role of parabolic subgroups in the theory of Artin groups, to explain why they matter, and to highlight key results obtained in recent decades.

\subsection{Some families of Artin groups}

As of now, there are very few techniques that apply to all Artin groups. As a result, group theorists work within different families, depending on their expertise. We enumerate some of these families below:

\begin{itemize}

\item \textbf{Right-Angled Artin Groups (RAAGs).} Also known as partially commutative groups, these are the Artin groups where the only possible relations between two generators are commutations.

\item \textbf{Spherical-type (or finite-type) Artin groups.} These are the Artin groups with finite associated Coxeter group. They are called spherical-type Artin group because these Coxeter groups are generated by reflections on the sphere. The irreducible spherical-type Artin group were classified by \cite{Coxeter} (\autoref{Dibujo:Spherical-type}).
\begin{figure}[h]
  \centering
  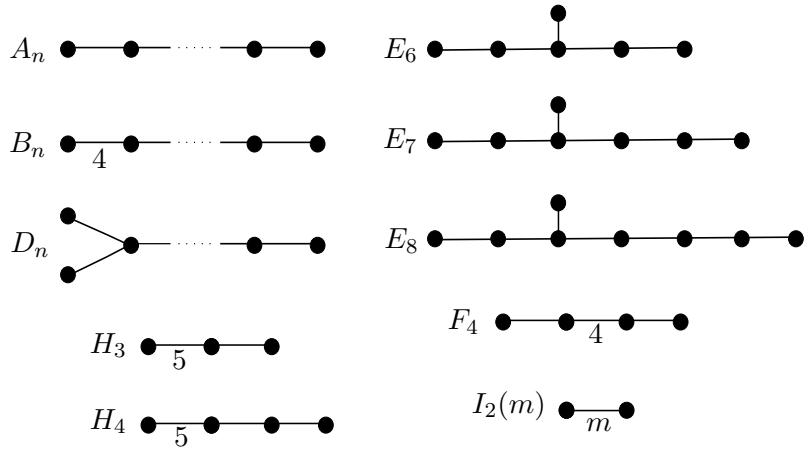
  \medskip
  \caption{The complete classification of irreducible spherical-type Artin group with the no-$2$-convention. The graph $A_n$ corresponds to the braid group on $n+1$ strands.}\label{Dibujo:Spherical-type}
\end{figure}

\item \textbf{FC-type Artin groups} (which include the two previously mentioned families) are characterized by the property that, under the no-$\infty$ convention, every clique in the Coxeter graph corresponds to a spherical-type Artin group. These groups can be viewed as amalgamated products of spherical-type Artin groups (see \autoref{exo:FC-decomposition}). They are precisely the Artin groups for which the Deligne complex is a flag complex \citep{CharneyDavis}, which is the origin of the name ``FC" (Flag Complex). For instance, the Artin group depicted in \autoref{Dibujo:convenciones} is of FC type.

\item \textbf{Affine-type (or Euclidean-type) Artin groups.} These are the Artin groups whose Coxeter groups correspond to reflection groups of Euclidean spaces (\autoref{Dibujo:Euclidean-type}).

\begin{figure}[h]
  \centering
  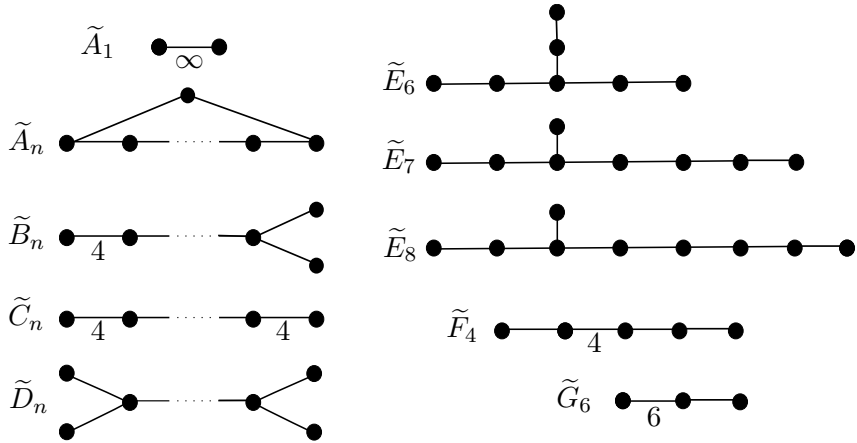
  \medskip
  \caption{The complete classification of irreducible Euclidean-type Artin groups with the no-$2$-convention. The graph $A_n$ corresponds to the Euclidean braid group.}\label{Dibujo:Euclidean-type}
\end{figure}

\item \textbf{Large-type Artin groups.} These are the Artin groups with no commutation relations.

\item \textbf{2-dimensional Artin groups} (including large-type Artin groups). These satisfy the condition that for every three generators $s,t,r\in S$, 
\[ \frac{1}{m_{s,t}}+\frac{1}{m_{s,r}}+\frac{1}{m_{t,r}}\leq 1, \]
which is equivalent to saying that the cohomological dimension of the group is at most 2. 

\end{itemize}

We can also extend this terminology to parabolic subgroups. For example, if \( P = \alpha A_X \alpha^{-1} \), where \( A_X \) is an Artin group of spherical type, we say that \( P \) is a parabolic subgroup of spherical type.

\section{Motivation}

\subsection{Parabolic subgroups and curves in the punctured disc}

\begin{figure}[h]
  \centering
    \def\svgwidth{0.7\columnwidth} 
  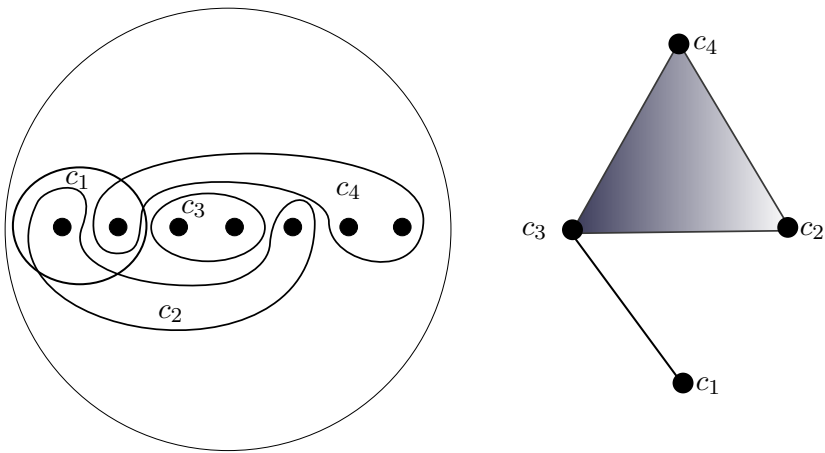
  \medskip
  \caption{A portion of the 1-skeleton of the curve complex of an \(n\)-punctured disc, corresponding to the four curves shown.}\label{Dibujo:CurveComplex}
\end{figure}

A reasonable question to explore in order to obtain results in Artin groups is: ``How many techniques from braid theory can be translated into Artin groups?" In particular, if we consider braids as the mapping class groups of the \(n\)-punctured disc \(\mathcal{D}_n\), we can extract a wealth of information from their action on the curve complex of \(\mathcal{D}_n\). The \emph{curve complex} is a flag complex. A flag complex is totally defined from its 1-skeleton: every set of $m$ vertices pairwise connected by an edge expand an $n$-dimensional simplex. The 1-skeleton of the curve complex of a surface is defined as follows: its vertices correspond bijectively to the isotopy classes of non-degenerate simple closed curves, and two vertices are adjacent if their corresponding classes admit disjoint representatives (see example in \autoref{Dibujo:CurveComplex}). It turns out that we can construct a analogous complex for every Artin group by using parabolic subgroups. 

\begin{figure}[h]
  \centering
    \def\svgwidth{0.7\columnwidth} 
\begingroup%
  \makeatletter%
  \providecommand\color[2][]{%
    \errmessage{(Inkscape) Color is used for the text in Inkscape, but the package 'color.sty' is not loaded}%
    \renewcommand\color[2][]{}%
  }%
  \providecommand\transparent[1]{%
    \errmessage{(Inkscape) Transparency is used (non-zero) for the text in Inkscape, but the package 'transparent.sty' is not loaded}%
    \renewcommand\transparent[1]{}%
  }%
  \providecommand\rotatebox[2]{#2}%
  \newcommand*\fsize{\dimexpr\f@size pt\relax}%
  \newcommand*\lineheight[1]{\fontsize{\fsize}{#1\fsize}\selectfont}%
  \ifx\svgwidth\undefined%
    \setlength{\unitlength}{547.65024741bp}%
    \ifx\svgscale\undefined%
      \relax%
    \else%
      \setlength{\unitlength}{\unitlength * \real{\svgscale}}%
    \fi%
  \else%
    \setlength{\unitlength}{\svgwidth}%
  \fi%
  \global\let\svgwidth\undefined%
  \global\let\svgscale\undefined%
  \makeatother%
  \begin{picture}(1,0.471566)%
    \lineheight{1}%
    \setlength\tabcolsep{0pt}%
    \put(0,0){\includegraphics[width=\unitlength,page=1]{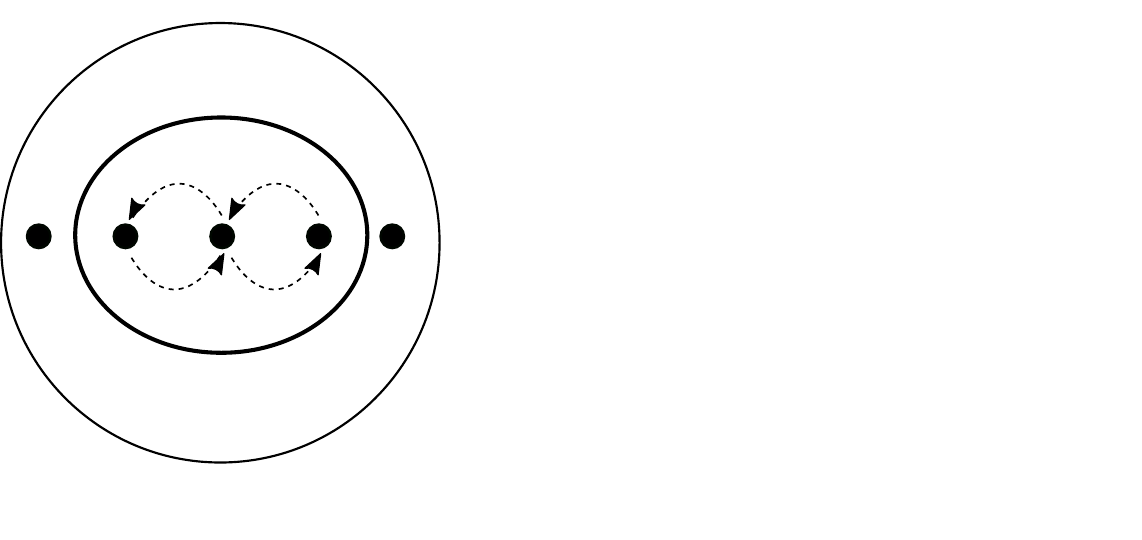}}%
    \put(0.13805763,0.25682765){\color[rgb]{0,0,0}\makebox(0,0)[lt]{\lineheight{1.25}\smash{\begin{tabular}[t]{l}$\sigma_2$\end{tabular}}}}%
    \put(0.22468839,0.25708707){\color[rgb]{0,0,0}\makebox(0,0)[lt]{\lineheight{1.25}\smash{\begin{tabular}[t]{l}$\sigma_3$\end{tabular}}}}%
    \put(0,0){\includegraphics[width=\unitlength,page=2]{paralelismo.pdf}}%
    \put(0.48136869,0.28038707){\color[rgb]{0,0,0}\makebox(0,0)[lt]{\lineheight{1.25}\smash{\begin{tabular}[t]{l}$\sigma_4^{-1}$\end{tabular}}}}%
    \put(0.16295792,0.02662048){\color[rgb]{0,0,0}\makebox(0,0)[lt]{\lineheight{1.25}\smash{\begin{tabular}[t]{l}$A_{\{\sigma_2,\sigma_3\}}$\end{tabular}}}}%
    \put(0.71325227,0.02602634){\color[rgb]{0,0,0}\makebox(0,0)[lt]{\lineheight{1.25}\smash{\begin{tabular}[t]{l}$\sigma_4^{-1}A_{\{\sigma_2,\sigma_3\}}\sigma_4$\end{tabular}}}}%
  \end{picture}%
\endgroup%

  \medskip
  \caption{In the braid group, irreducible parabolic subgroups are in bijection with the set of curves up to isotopy in the punctured disc.}\label{Dibujo:paralelismo}
\end{figure}

We say that a (class of) curve(s) is \emph{standard} if there is a representative that intersects the diameter of the disc ---the one containing all the punctures of $\mathcal{D}_n$--- only at two points.
First, note that the vertices of the curve complex of $\mathcal{D}_n$ are in bijection with the irreducible parabolic subgroups of $\mathcal{B}_n$. In more detail, each irreducible standard parabolic subgroup $A_X$ can be associated with a unique standard curve $c$ that encloses the punctures involved in the generators of $X$ (and vice versa). Moreover, every irreducible parabolic subgroup $\alpha^{-1}A_X \alpha$ is uniquely associated with the image of $c$ under the action of $\alpha$ (and the converse also holds, because any curve can be viewed as a standard curve deformed by a braid). (See \autoref{Dibujo:paralelismo} for an illustration.)

The \emph{complex of irreducible parabolic subgroups} \citep{CGGW} for an Artin group is a flag complex whose vertices are irreducible parabolic subgroups. Two vertices, representing subgroups $P$ and $Q$, are adjacent if and only if one of the following three conditions holds:
\[
P \subset Q, 
\quad
Q \subset P,
\quad
\text{or}
\quad
\bigl(P \cap Q \text{ is trivial and } pq=qp \text{ for every }(p,q)\in P\times Q\bigr).
\]
Checking this adjacency condition can be challenging. The following results make it simpler to handle in the spherical cases, which are the only ones conjectured to have non-trivial center:

\begin{theorem}[\citealp{CGGW,Rose}]\label{thm:ccgw-mr}
For FC-type Artin groups, the previous adjacency condition is equivalent to the equality \( z_Q z_P = z_P z_Q\), where \(z_P\) and \(z_Q\) are the generators of the centres of \(P\) and \(Q\) (see \autoref{section:Garside}), respectively, and \(P, Q\) have spherical type.
\end{theorem}

\begin{question}\label{q:adjacency-condition}
Is there any effective way to check the adjacency condition in general?
\end{question}

We know that the curve complex of \(\mathcal{D}_n\) is Gromov-hyperbolic, but what do we know for other cases?

\begin{theorem}[\citealp{CalvezCisneros}]\label{thm:calvez}
The complex of irreducible parabolic subgroups is Gromov-hyperbolic for spherical-type \(B\) and Euclidean types \(\widetilde{A}\) and \(\widetilde{C}\).
\end{theorem}

\begin{question}\label{q:gromov-hyperbolic}
In which cases is the complex of irreducible parabolic subgroups Gromov-hyperbolic?
\end{question}

Also, for mapping class groups we have the \emph{Nielsen--Thurston classification}, which effectively categorizes braids into three types:
\begin{itemize}
    \item \textbf{Periodic}, if a power of the braid is central;
    \item \textbf{Reducible}, if a power of the braid preserves a non-trivial family of curves;
    \item \textbf{Pseudo-Anosov}, if no power of the braid preserves any curve in \(D_n\).
\end{itemize}

This classification, together with the notion of a canonical reduction system, is a powerful tool —for a classic reference, see \citep{FarbMargalit}—. In the reducible case, one can pick a family of curves and examine the different regions of the surface delimited by these curves. The element will act in these parts either periodically or in a pseudo-Anosov manner. The pseudo-Anosov case has a far richer description: in this scenario, there are two transverse measurable foliations that are respectively stretched or widened by a factor related to the braid, providing substantial structure. Since periodic braids can be considered a trivial case, the canonical reduction system allows us to extract information from reducible braids by studying the other two cases. Observe that we can devise a similar classification for parabolic subgroups.

\begin{definition}[Nielsen-Thurston Classification for Artin groups]
An element $g$ of an Artin group can be 

\begin{itemize}
    \item \textbf{Periodic}, if a power of $g$ is central;
    \item \textbf{Reducible}, if a power of $g$ preserves a non-trivial family of parabolic subgroups under conjugacy;
    \item \textbf{Pseudo-Anosov}, if no power of the braid preserves any  parabolic subgroup under conjugacy.
\end{itemize}

\end{definition}

\begin{question}
Is there a notion of canonical reduction system for Artin groups? 
\end{question}

The questions above are very hard, and depend on understanding more basic properties of parabolic subgroups.

\subsection{Parabolic subgroups and the $K(\pi,1)$ conjecture}

One of the most important long-standing conjectures in the theory of Artin groups is the $K(\pi,1)$ conjecture, which establishes a powerful connection between algebra and topology. 
Roughly speaking, every Artin group is known to be the fundamental group of the complement $M$ of a certain hyperplane arrangement. The conjecture asserts that $M$ is a Eilenberg–MacLane space of type $K(\pi,1)$ for the Artin group; that is, $\pi_m(M) = 0$ for all $m > 1$. Equivalently, this means that the universal cover of $M$ is contractible.
For a detailed survey of this conjecture, see \citep{ParisK} and the more recent survey by Rachael \cite{Boyd}.

\medskip
To begin with a simple example, consider the pure braid group. Viewed from the perspective of mapping class groups, pure braids describe motions of the punctures in \(\mathcal{D}_n\) such that each puncture returns to its original position. Since the motion is a braid, no two punctures may occupy the same point at the same time. This interpretation allows us to realize the pure braid group on \(n\) strands as the fundamental group of the configuration space:
\[
\pi_1\left(\mathbb{C}^n \setminus \bigcup_{i \ne j} \{z_i = z_j\}\right).
\]
For the full braid group, where punctures may start and end at different positions, we quotient this space by the symmetric group \(\mathcal{S}_n\), obtaining:
\[
\mathcal{B}_n = \pi_1\left(\frac{\mathbb{C}^n \setminus \bigcup_{i \ne j} \{z_i = z_j\}}{\mathcal{S}_n}\right).
\]
It was shown by \cite{Fox} that this space is a \(K(\pi, 1)\) space for the braid group.
In general, it was proven by \cite{Vanderlek} that, given an Artin group \(A\), its associated Coxeter group \(W\) acts faithfully on a non-empty open convex cone \(I\), such that the union of regular orbits is the complement in \(I\) of a finite collection \(\mathcal{H}\) of linear hyperplanes. Moreover, the Artin group \(A\) is the fundamental group of the space:
\[
M = \frac{(I \times I) \setminus \bigcup_{H \in \mathcal{H}} H \times H}{W}.
\]

The progress toward proving that \( M \) is an Eilenberg--MacLane space has occurred at discrete points in time. \cite{Brieskorn} (partially) and \cite{Deligne} proved the conjecture for spherical-type Artin groups. \cite{Okonek} proved it for the affine types \( \tilde{A} \) and \( \tilde{C} \). \cite{Hendriks} established the result for large-type Artin groups, and this was later extended to the 2-dimensional case by \cite{CharneyDavis2}. \cite{CharneyDavis} also generalized the spherical-type case to FC-type Artin groups. After that, there was a gap until 2010, when \cite{CMS} proved the conjecture for the affine type \( \tilde{B} \). Building on the work of \cite{McCammondSulway}, \cite{PaoliniSalvetti} finally proved the conjecture for all affine Artin groups.
In a significant part of these proofs, a model homotopy equivalent to \( M \) is constructed. The main model used is the \emph{Salvetti complex}, and to model the universal cover we use the \emph{Deligne complex}. Parabolic subgroups play a role in the construction of this complex. Also, \cite{Ellis} proved that in an Artin group all the standard parabolic subgroups without $\infty$-relations satisfy the $K(\pi,1)$ conjecture, then the whole group also fulfils it. 

\medskip

\begin{figure}[h]
  \centering
    \def\svgwidth{0.6\columnwidth} 
  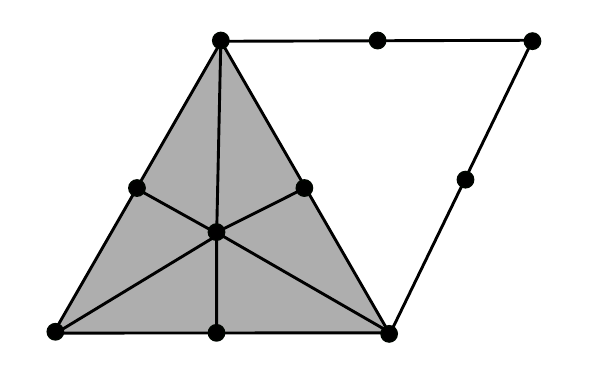
  \medskip
  \caption{The derived complex of the set  $$\{\{a\},\{b\},\{c\},\{d\},\{a,b\},\{a,c\}, \{a,d\}, \{b,c\},\{c,d\},\{a,b,c\}\}$$ partially ordered by inclusion.}\label{Dibujo:derived}
\end{figure}

To understand the construction, we first need to know what is the derived complex associated with a partially ordered set (poset). Given a poset, associate to every chain of inequalities of $k$ elements a simplex of dimension $k$. (See \autoref{Dibujo:derived} for an example.)

Now, consider an Artin group \( A \) with associated Coxeter group \( W \), and let \( S^f \) be the set consisting of all the subsets of generators of \( W \) that generate a finite Coxeter subgroup (which is analogous to taking the subset of generators of \( A \) that generate a standard parabolic subgroup of spherical type). We can define a partial order \( \preceq \) on \( W \times S^f \) as follows: \( (u, X) \preceq (v, Y) \) if \( X \subset Y \), \( v^{-1}u \in W_Y \), and \( v^{-1}u \) is minimal in the coset \( v^{-1}uW_X \). The derived complex of this poset is what we call the \emph{Salvetti complex} associated with \( A \). The reader can find more details about this complex in \autoref{section:retractionSalvetii}.
Next, we can partially order by inclusion the set 
\[
AS^f = \{ \alpha A_T \mid \alpha \in A,\, T \in S^f \},
\] 
and the corresponding derived complex is the \emph{Deligne complex}.  
For example, the complex in \autoref{Dibujo:derived} could correspond to the portion of the Deligne complex associated with the cosets \( 1 \cdot S^f \) of the Artin group
\[
\langle a,b,c,d \mid aba=bab,\; aca=cac,\; adad=dada,\; bc=cb,\; cdc=dcd \rangle.
\]
The Deligne complex can also be constructed as a complex of groups, where the fundamental domain is built from the inclusions of spherical-type parabolic subgroups, and the entire complex is obtained via the action of the group (see \citep[Chapter~III.$\mathcal{C}$]{BridsonHaefliger} for the construction).  
Under this construction, the spherical-type parabolic subgroups of \( A \) are precisely the stabilizers of the simplices of the Deligne complex.  

If, instead of considering only spherical-type parabolic subgroups, we take all non-trivial parabolic subgroups, we obtain in an analogous way a complex in which every parabolic subgroup stabilizes some simplex. This complex is called the \emph{Artin complex}, and it is precisely defined as the complex whose barycentric subdivision is the derived complex of the set   
\[
A\mathcal{S}=\{ \alpha A_T \mid \alpha \in A,\, T \text{ is a non-trivial parabolic subgroup of } A \}.
\] 
In the Artin complex, the portion corresponding to \( 1 \cdot \mathcal{S} \) gives a \( (n-1) \)-dimensional simplex, where $n$ is the number of standard generators. For example, for an Artin group with three generators $\{a,b,c\}$ we obtain a 2-dimensional simplex (the shadowed simplex that we can see in \autoref{Dibujo:derived} and in \autoref{Dibujo:ArtinComplex}). This is in fact the fundamental domain of the complex for the action of the group. We will explain this further in \autoref{section:GGT}.

\section{Tools to work with parabolic subgroups}

In this section we want to show some different approaches to work in Artin groups, explaining how these tools have been used to achieve different important results about parabolic subgroups. In the first subsection we will deal with a basic question:

\begin{conjecture}\label{Conjecture:Intersection}
The intersection of two parabolic subgroups is again a parabolic subgroup. 
\end{conjecture}

In \autoref{section:Garside}, we will present the proof for spherical-type Artin groups \citep{CGGW}, which was the first major breakthrough on the problem beyond the case of RAAGs \citep{Duncan}. In \autoref{section:GGT}, we will explore a geometric group theory approach using the Artin and Deligne complexes, where parabolic subgroups are interpreted as stabilizers of simplices. This viewpoint leads to proofs of the conjecture for spherical-type parabolic subgroups in FC-type Artin groups \citep{Rose}, for large-type Artin groups \citep{CMV}, and for certain 2-dimensional Artin groups \citep{Blufstein}. We will also see how Bass--Serre theory can be applied to prove the result in FC-type Artin groups when one of the parabolic subgroups is of spherical type \citep{Moller}.
In \autoref{section:retractions}, we will study the use of algebraic retractions, which have been employed to prove the conjecture for even FC-type Artin groups \citep{AntolinFoniqi}. For Euclidean-type Artin groups, we currently have proofs for types $\widetilde{A}$ and $\widetilde{C}$, due to the work of \cite{Haettel}. In \autoref{section:euclidean}, we will revisit the case of Euclidean braid groups, following the argument in \citep{CGP}.

\medskip

It is not difficult to take the previous conjecture one step further: 

\begin{proposition}
If $A_S$ is an Artin group satisfying \autoref{Conjecture:Intersection}, then the intersection of an arbitrary number of parabolic subgroups is a parabolic subgroup. 
\end{proposition}

\begin{proof}

The proof follows the structure of \citep[Proposition 10.1]{CGGW} and \citep[Corollary~16]{CMV}.
Let $\mathcal{P}$ be an arbitrary collection of parabolic subgroups of $A_S$, and set \( Q = \bigcap_{P \in \mathcal{P}} P \). Since \( Q \) is contained in every parabolic subgroup in \( \mathcal{P} \), in order to apply \autoref{Conjecture:Intersection}, it suffices to show that \( Q \) coincides with a finite intersection of parabolic subgroups.

Observe that every parabolic subgroup is a conjugate of a standard parabolic subgroup. As~\( A_S \) is countable and the set of standard parabolic subgroups is finite, the set of all parabolic subgroups of \( A_S \) is also countable. In particular, we may assume that \( \mathcal{P} = \{P_1, P_2, P_3, \dots\} \) is countable. Define
\[
Q_m = \bigcap_{1 \leq i \leq m} P_i.
\]
If \autoref{Conjecture:Intersection} holds, then each \( Q_m \) is a parabolic subgroup. Since \( Q = \bigcap_{i \in \mathbb{N}} Q_i \), it is enough to prove that the set \( \{Q_m \mid m \in \mathbb{N}\} \) is finite.

Note that this defines a descending chain
\[
Q_1 \supseteq Q_2 \supseteq Q_3 \supseteq \dots
\]
By \autoref{theorem_restandardisation} below, after conjugating so that each \( Q_i \) becomes a standard parabolic subgroup, we see that \( Q_i \) (for \( i > 1 \)) is a parabolic subgroup of \( Q_{i-1} \). This implies that the rank of~\( Q_i \) is strictly less than that of \( Q_{i-1} \). Since \( |S| \) is finite, this strict inequality cannot continue indefinitely, and the chain must eventually stabilize. This proves the claim.
\end{proof}

Finally, in \autoref{section:retractionSalvetii}, we will see how the Salvetti complex can be used to prove two of the few general results known for Artin groups. The first one states that if we take a geodesic word (with respect to the word length in the standard generators of $A_S$) representing an element of a parabolic subgroup \( A_X \subseteq A_S \), then all the letters in the word must lie in \( X \sqcup X^{-1} \).

\begin{theorem}[\citealp{CharneyParis}]\label{theorem_charneyparis}
Standard parabolic subgroups are convex with respect to the word length. 
\end{theorem}

The second theorem roughly says that a parabolic subgroup $P$ contained in a standard parabolic subgroup $A_X$ is a parabolic subgroup of $A_X$.

\begin{theorem}[\citealp{BlufsteinParis}]\label{theorem_restandardisation}
Let $P$ and $A_X$ be two parabolic subgroups of an Artin group $A$ such that $P\subset A_X$. Then there are $\alpha\in A_X$ and $Y\subset X$ such that $P=\alpha A_Y \alpha^{-1}$. 
\end{theorem}

This result had been previously proved for spherical-type Artin groups \citep{Godelle,Paris} and extended to FC-type Artin groups \citep{GodelleFC}. Also it had been proven for 2-dimensional \citep{Godelle2dim} and Euclidean types $\widetilde{A}$ and $\widetilde{C}$ \citep{Haettel}.

\subsection{Intersection of parabolic subgroups}

\subsubsection{The spherical-type case: Garside Theory}\label{section:Garside}

Spherical-type Artin groups have a very nice underlying structure called Garside structure. It has its origin in the solution to the word problem for braids \citep{Adyan,Garside,Elrifai} and it was first used for all spherical-type Artin groups in \citep{BrieskornSaito}. Later, \cite{DehornoyParis} isolated these properties into what we call now Garside groups. 

\medskip

For a spherical-type Artin group $A$ with standard set of generators $S$ we need to consider the following elements and properties: 

\begin{itemize}

\item The monoid $A^+$ of positive elements of $A$. 

\item A partial order $\po$ defined as follows: $a\po b$ if and only if there is $c\in A^+$ such that $b=ac$. This order is called \emph{prefix order} and we say that $a$ is a prefix of $b$. For all $a,b\in A$ there is a unique least common multiple $a\vee b$ and a unique great common divisor $a \wedge b$. This order is invariant under left multiplication.

\item The Garside element $\Delta := \bigvee_{s\in S} s$, i.e. the least common multiple of the collection of standard generators of $A$, satisfies $\Delta A^+ \Delta^{-1}= A^+$. 

\item The prefixes of $\Delta$ are called \emph{simple elements} and also generate the group.

\item $A$ is atomic: The length of words written in $S$ (called also \emph{atoms}) representing a positive element has an upper bound (it is actually always the same length).

\end{itemize}

One can also define a partial \emph{suffix order} $\soo$ such that $a\soo b$ if and only if there is $c\in A^+$ such that $a=cb$ with the same properties. In this case we say that $b$ is a suffix of $a$. We will also have that $\Delta$ is the least common multiple of the atoms with respect this order. By default, we will work with the prefix order. 

\medskip

A standard parabolic subgroup $A_X$ is also a spherical-type Artin group, so it has its own Garside structure, with Garside element $\Delta_X$. It is well-known that for irreducible spherical-type Artin group either the Garside element or its square generates the center, we will call $z_X$ the generator of the center of $A_X$. It is easy to see that for every parabolic subgroup $P=\alpha^{-1}A_X\alpha$, $z_P:=\alpha^{-1}z_X \alpha$ is the generator of $Z(P)$.

\medskip

There is a quick way to compute the Garside element.  
By \citep{BrieskornSaito}, we know that the power of the Garside element that generates the center is always a power of the product of the standard generators (in any order).  
For example, for $A_2$ we have 
$$
\Delta^2 = (s_1 s_2)^{3} = s_1s_2s_1s_2s_1s_2,
$$
where $\Delta=s_1s_2s_1$ is the least common multiple of $s_1$ and $s_2$.
In general, if $A$ has $n$ generators, the generator of its center will be $(s_1s_2\cdots s_n)^{k_\Delta}$. The exponent~\( k_\Delta \) depends on the so-called Coxeter number, and these numbers are well known in the spherical case (see \autoref{tab:kdelta} and Section~3.18 in \citealp{Humphreys} for a reference).  
We refer the reader to \citep{ParisTabla} for further details.  
However, it is also possible to compute \(\Delta\) directly by hand, using 
$
\Delta = \bigvee_{s \in S} s.$ We refer the reader to solution of \autoref{exo:GarsideElement} to see explicit computations of least common multiples.

\begin{table}[h!]
\centering
\begin{tabular}{llc}
\toprule
\(\Gamma\) & Conditions & \(k_\Delta\) \\
\midrule
\(A_1\) &  & \(1\) \\
\(A_n\) & \(n \ge 2\) & \(n+1\) \\
\(B_n\) & \(n \ge 2\) & \(n\) \\
\(D_n\) & \(n \ge 4\) even & \(n-1\) \\
\(D_n\) & \(n \ge 5\) odd & \(2n-2\) \\
\(E_6\) &  & \(12\) \\
\bottomrule
\end{tabular}
\hspace{1.5cm} 
\begin{tabular}{llc}
\toprule
\(\Gamma\) & Conditions & \(k_\Delta\) \\
\midrule
\(E_7\) &  & \(9\) \\
\(E_8\) &  & \(15\) \\
\(F_4\) &  & \(6\) \\
\(H_3\) &  & \(5\) \\
\(H_4\) &  & \(15\) \\
\(I_2(p)\) & \(p \ge 6\) even & \(p/2\) \\
\(I_2(p)\) & \(p \ge 5\) odd & \(p\) \\
\bottomrule
\end{tabular}
\caption{Values of \(k_\Delta\) for irreducible spherical-type Coxeter systems}
\label{tab:kdelta}
\end{table}

\begin{exercise}\label{exo:GarsideElement}
Compute the Garside element of $A[H_3]$. \emph{(See solution in page \pageref{sol:Garside}.)}
\end{exercise}

Thanks to the Garside structure, any element \(\alpha \in A\) has normal form:

\begin{theorem}[\citealp{Elrifai,Thurston1992,Adyan}]
Any element $\alpha$ of a Garside group can be written as
\[
\alpha = \Delta^p x_1 \cdots x_r,
\]
where each \(x_i\) is a simple element and \(x_i x_{i+1} \wedge \Delta = x_i\) (we say that $x_i\cdot x_{i+1}$ is \emph{left-weighted}).  
This is called the \emph{left normal form} of~\(\alpha\).  

Analogously, using the suffix order, we can write \(\alpha\) as
\[
\alpha = x'_r \cdots x'_1 \Delta^p,
\]
where each \(x'_i\) is simple and \(x'_{i+1} x'_i \wedge \Delta = x'_i\).  
This is called the \emph{right normal form} of \(\alpha\).
\end{theorem}
 
We define the \emph{infimum} of \(\alpha\) as \(\mathrm{inf}(\alpha) := p\),  
the \emph{supremum} of \(\alpha\) as \(\mathrm{sup}(\alpha) := p + r\),  
and the \emph{length} of \(\alpha\) as \(\ell(\alpha)\).  
It is known that these numbers are independent of whether the left or right normal form is used  
(see Remarque~1.10 in \citealp{CumplidoThese} for an explanation).

\medskip

Sometimes it is convenient to express this normal form in a slightly different way,  
using Charney’s \emph{mixed normal form} \citep{Charney95}.  
If \(\alpha\) is entirely positive or entirely negative, we simply write \(\alpha\).  
Otherwise, set
\[
a^{-1} = \Delta^p x_1 \cdots x_{-p}
\quad\text{and}\quad
b = x_{-p} \cdots x_r,
\]
and write \(\alpha = a^{-1} b\).  
This is called the \emph{np-normal form} (negative–positive) and has the property that \(a \wedge b = 1\).  
We can similarly use the suffix order to obtain \(\alpha = a b^{-1}\),  
the \emph{pn-normal form} (positive–negative).

\begin{proposition}[{\citealp[Lemma 4.5]{Elrifai}}]\label{Prop:negativo} 
If $\alpha=\Delta^p x_1\cdots x_r$ is in left normal form then, if we set $x_i'=\Delta x_i^{-1}$
$$\alpha^{-1}= \Delta^{-(p+r)}\cdot \left( \Delta^{p+r-1} x'_r \Delta^{-(p+r-1)}\right)\cdot \left( \Delta^{p+r-2} x_{r-1}' \Delta^{-(p+r-2)}\right)\cdots\left ( \Delta^{p} x'_{1} \Delta^{-p}\right)$$ is also in left normal form.  
\end{proposition}
 
\begin{exercise}\label{exo_inf_sum_mixednoralforms}
Let $\alpha=a^{-1}b$ be in mixed np-normal form. Prove that, if $p<0$, $\mathrm{sup}(a)=-\mathrm{inf}(\alpha)$ and $\mathrm{sup}(b)=\mathrm{sup}(\alpha)$. \emph{(See solution in page \pageref{sol:inf}.)}
\end{exercise}

Our aim will be to provide the necessary ingredients to understand the proof of the following theorem.  
Although the proof is long and technical, our goal is to give the reader a clear picture of its overall structure  
and of how the key technical results are used throughout.

\begin{theorem}[{\citealp[Theorem 9.5]{CGGW}}]\label{theorem:intersection}
Let $P$ and $Q$ be two parabolic subgroups of a spherical-type Artin group $A$. Then $P\cap Q$ is a parabolic subgroup.
\end{theorem}

The proof is based on a weaker result:

\begin{theorem}[{\citealp[Theorem 1.1]{CGGW}}]\label{theorem:unique_parabolic_closure}
Every element $x$ of a spherical-type Artin group $A$ has a parabolic closure $P_x$. That is, every $x$ is contained in a unique minimal parabolic subgroup with respect to the inclusion.
\end{theorem}

\begin{proof}
The proof we present here follows the approach in \citep{CGGW}, but we explain a theoretically simpler version based on the swap operation introduced in \citep{MenesesMarin} to generalize the argument to certain Garside groups. We give here the skeleton of the proof.

We start by defining an operation $\phi$, called a \emph{swap}, on any element $x$ written in its mixed normal form $x = a^{-1} b$ as
\[
\phi(x) := b a^{-1}.
\]
It is shown in \citep[Proposition~4.13]{MenesesMarin} that for every $x$ there exist positive integers $m < n$ such that $\phi^m(x) = \phi^n(x)$. The set
\[
\{\phi^m(x), \phi^{m+1}(x), \dots, \phi^{n-1}(x)\}
\]
is called a \emph{circuit for the swap}. Any element in such a circuit is called a \emph{recurrent element}.  

Let $\mathrm{R}(x)$ denote the set of all recurrent elements conjugate to $x$. By \citep[Theorem~4.29]{MenesesMarin}, each $y \in \mathrm{R}(x)$ has a \emph{parabolic closure}
\[
P_y := A_{\mathrm{Supp}(y)},
\]
where $\mathrm{Supp}(y)$ is the set of all standard generators that appear in the mixed normal form $c^{-1}d$ of $y$ (see the sketch below). This support is well defined because, in Artin groups, the relations are homogeneous so a positive element is always represented by positive words involving the same generators.

Now suppose $x = \alpha^{-1} y \alpha$ with $y \in \mathrm{R}(x)$. We want to show that
\[
P_x = \alpha^{-1} P_y \alpha.
\]
Indeed, if $Q$ is any parabolic subgroup containing $x$, then $y \in \alpha Q \alpha^{-1}$, hence $P_y \subseteq \alpha Q \alpha^{-1}$, and conjugating back gives $P_x \subseteq Q$.

\medskip
\noindent\emph{Skeleton of the proof of \citep[Theorem~4.29]{MenesesMarin}.}
Let $y$ be a recurrent element. We claim that $P_y = A_{\mathrm{Supp}(y)}$. Since $A_{\mathrm{Supp}(y)}$ is already a standard parabolic subgroup, it remains to prove its minimality.

Let $Q = \alpha^{-1} A_X \alpha$ be a parabolic subgroup containing $y$. Then $\alpha y \alpha^{-1} \in A_X$. By repeated applications of the swap operation, we obtain a recurrent element $\beta \alpha y \alpha^{-1} \beta^{-1}$ for some $\beta \in A$. We may assume $\beta \in A_X$, because $A_X$ is also a spherical-type Artin group, and the mixed normal forms in $A$ and $A_X$ coincide (see \citealp[Section~3]{CGGW} for an explanation).

Thus $\beta \alpha y \alpha^{-1} \beta^{-1} \in A_X$, and all the letters in its mixed normal form correspond to generators in $X$. Let
\[
Y := \mathrm{Supp}\big( \beta \alpha y \alpha^{-1} \beta^{-1} \big) \subseteq X.
\]
By \citep[Proposition~4.28]{MenesesMarin}, if $y$ and $z$ are recurrent elements and $\gamma$ is such that $\gamma y \gamma^{-1} = z$, then
\[
\gamma \, A_{\mathrm{Supp}(y)} \, \gamma^{-1} = A_{\mathrm{Supp}(z)}.
\]
Applying this with $z = \beta \alpha y \alpha^{-1} \beta^{-1}$ yields
\[
A_Y = \beta \alpha \, A_{\mathrm{Supp}(y)} \, \alpha^{-1} \beta^{-1} \subseteq A_X.
\]
Conjugating back, we obtain
\[
A_{\mathrm{Supp}(y)} \subseteq \alpha^{-1} \beta^{-1} A_X \beta \alpha
= \alpha^{-1} A_X \alpha = Q,
\]
which proves the minimality of $A_{\mathrm{Supp}(y)}$.

\end{proof}

\begin{exercise}\label{exo_parabolicclosure}
Prove that for every parabolic subgroup $Q$, we have that $P_{z_Q}=Q$. \emph{(See solution in page \pageref{sol:parabolicclosure}.)}
\end{exercise}

\begin{exercise}\label{exo_recurrent_swapping}
Prove that when we apply recurrent swapping to an element, the infimum can only increase and the supremum can only decrease. \emph{(See solution in page \pageref{sol:recurrent}.)}
\end{exercise}
\begin{lemma}\label{exo_birman_ko_lee} 
Let \(\alpha\) be an element in a spherical-type Artin group.  
If \(\alpha\) is conjugate to a positive element, then we can reach a positive element by conjugating by an element \(c'\) such that 
\[
\mathrm{sup}(c') \leq |\mathrm{inf}(\alpha)| \, l(\Delta).
\] 
\end{lemma}

\begin{proof}
If \(\alpha\) is already positive, the result is immediate.  
Suppose instead that \(\alpha\) is not positive, so that its infimum is not maximal.  

Given the left normal form \(\alpha = \Delta^p x_1 \cdots x_r\), we define the \emph{cycling} of \(\alpha\) as
\[
c(\alpha) = \Delta^p x_2 \cdots x_r \, (\Delta^p x_1 \Delta^{-p}) 
= (\Delta^p x_1^{-1} \Delta^{-p}) \, \alpha \, (\Delta^p x_1 \Delta^{-p}).
\]
Note that $\Delta^p x_1 \Delta^{-p}$ is either $x_1$ or $\Delta x_1 \Delta^{-1}$ and it is a simple element, so the cycling cannot increase the length of the normal form and it cannot decrease the infimum.

If \(\mathrm{inf}(\alpha)\) is not maximal among the infimums of the elements in the conjugacy class of $\alpha$, then by \citep[Theorem~1]{BirmanKoLee2} the infimum can be increased after at most \(l(\Delta)-1\) cycling operations. Since we need to increase the infimum by \(|\mathrm{inf}(\alpha)|\), the conjugating element \(c'\) is a product of at most \(|\mathrm{inf}(\alpha)| \, l(\Delta)\) simple elements, as required.  
\end{proof}

\begin{proof}[Proof of \autoref{theorem:intersection}] Given an element \(g\) in an Artin group, we denote by \(l(g)\) its word length in the standard generators of the presentation.  
We may assume that \(P\) and \(Q\) are neither trivial nor the whole group, and that their intersection is non-trivial, so there exists some \(x \in P \cap Q\).  
We choose \(x\) so that its parabolic closure \(P_x \subseteq P \cap Q\) is exactly \(P \cap Q\).  
To achieve this, we impose a maximality condition. If \(P_x = \alpha^{-1} A_Y \alpha\), define
\[
\varphi(x) := l(\Delta_Y),
\]
which is well defined by \citep[Proposition~9.4]{CGGW}.  
We then select \(x\) with maximal \(\varphi(x) =: n\) among all elements of \(P \cap Q\).  
Up to conjugacy, we may assume that \(P_x = A_Z\) is standard.  
In particular, \(\Delta_Z \in P \cap Q\).

We aim to show that any \(w \in P \cap Q\) lies in \(A_Z\).  
It suffices to prove that the parabolic closure \(P_w := T\) of \(w\) is contained in \(A_Z\).  
Since \(P_{z_T} = T\) (\autoref{exo_parabolicclosure}), it is enough to show that \(z_T \in A_Z\).  
For this purpose, we consider the family of elements
\[
\beta_m := z_T \, \Delta_Z^m \in P \cap Q.
\]
Our goal, to complete the proof, is to show that for \(m\) large enough, \(\beta_m \in A_Z\).

\medskip

\noindent
\emph{Claim.}
For $m>M$ big enough, we always have that $\beta_m$ is conjugate to a positive element $\widetilde{\beta}_m$.

\noindent
\emph{Proof of the claim.}
Consider the element \(z_T = a^{-1} b\) in its mixed normal form.  
If we write \(a = a_1 \cdots a_r\) and \(b = b_1 \cdots b_s\) in left normal form, then
\[
\beta_m = a_r^{-1} \cdots a_1^{-1} \, b_1 \cdots b_s \, \Delta_Z^m.
\]
Observe that the mixed normal form of \(\beta_m\) is obtained by cancelling, in the middle, the largest common prefix between \(a\) and \(b \Delta_Z^m\).  
This implies that \(\mathrm{inf}(\beta_m) \geq -r\) and \(\mathrm{sup}(\beta_m) \leq s + m\).
Using recurrent swappings, we can conjugate \(\beta_m\) to a recurrent element \(\widetilde{\beta}_m\).  
By \autoref{exo_recurrent_swapping}, recurrent swapping can only increase the infimum and decrease the supremum, hence \(\mathrm{inf}(\widetilde{\beta}_m) \geq -r\) and \(\mathrm{sup}(\widetilde{\beta}_m) \leq s + m\).
By \autoref{exo_inf_sum_mixednoralforms}, this means that if \(\widetilde{\beta}_m = x_m^{-1} y_m\) is in mixed normal form, then \(x_m\) has at most \(r\) factors in its normal form decomposition and \(y_m\) has at most \(s + m\) factors.  
We want to prove that \(x_m = 1\) for sufficiently large \(m\).

Let \(U_m := \mathrm{Supp}(\widetilde{\beta}_m)\).  
By the maximality condition on \(\varphi(x)\), we have \(|\Delta_{U_m}| \leq n\).  
The normal form decomposition of \(\widetilde{\beta}_m\) coincides with its normal form decomposition when considering only the Garside structure of \(A_{U_m}\) (see \citealp[Section~3]{CGGW} for details).  
It follows that each factor in the normal form of \(\widetilde{\beta}_m\) has length at most \(n\), and all factors in \(x_m\) and \(y_m\) belong to \(A_{U_m}^+\).
We will show that, for sufficiently large \(m\), the element \(y_m\) contains \(\Delta_{U_m}\) as a factor, making it the first factor in its left normal form.  
Since there can be no cancellations between \(x_m^{-1}\) and \(y_m\), this will imply \(x_m = 1\), as desired.

Suppose that \(y_m\) does not contain any \(\Delta_{U_m}\).  
In that case, each of its factors has length at most \(n-1\).  
Let \(\xi(g)\) denote the \emph{exponent sum} of an element \(g\) in the standard generators,  
that is, the sum of the exponents of all letters in any word representing \(g\) in the Artin presentation.  
This is well defined in Artin groups and invariant under conjugacy.
In particular, \(\xi(\beta_m) = \xi(\widetilde{\beta}_m)\).
On the one hand, by construction of \(\beta_m = a^{-1} b \, \Delta_Z^m\), we have
\[
\xi(\beta_m) = l(b) - l(a) + n m = n m + K,
\]
where \(K := l(b) - l(a)\) is a constant.  
On the other hand, if \(y_m\) contains no \(\Delta_{U_m}\), then each of its factors has length at most \(n-1\),  
and thus
\[
\xi(\widetilde{\beta}_m) = l(y_m) - l(x_m) \le l(y_m) \le (n-1)(s + m) = (n-1)m + k,
\]
where \(k\) is a constant.
Comparing the two expressions for the exponent sum gives
\[
n m + K \le (n-1)m + k,
\]
that is,
\[
m \le k - K.
\]
Therefore, for every \(m > k - K=:M\), the element \(y_m\) must contain \(\Delta_{U_m}\) as a factor.  
This implies that \(\Delta_{U_m}\) is the first factor in the left normal form of \(y_m\), as required.

\medskip
 
Suppose that \( m \) is large enough so that \( \beta_m \) is always a conjugate of a positive element.  According to \autoref{exo_birman_ko_lee}, we can select  \( c_m \in A_S\) so that
\[
\widetilde{\beta}_m := c_m^{-1} \beta_m c_m
\]
 is positive and 
\[\mathrm{sup}(c') \leq  |\mathrm{inf}(\beta_m)| \cdot l(\Delta_S) = r \, l(\Delta_S),
\]
where \( r \) is independent of \( m \).  
This imply that the set of values taken by the sequence \( \{ c_i \}_{i \ge 1} \) is finite.  

\medskip

Next, we prove that if \( R_m = s_1 \cdots s_{N_m} \) denotes the non-\(\Delta_{U_m}\) part of the normal form of \(\widetilde{\beta}_m\) inside \(A_{U_m}\), then the sequence \(\{R_i\}_{i \geq 1}\) is finite.  
It suffices to show that \(N_m\) is bounded above by a constant independent of \(m\).
From the previous claim, \(\widetilde{\beta}_m\) has at most \(s+m\) factors in its normal form, where \(s\) is the supremum of the normal form of \(b\) in the mixed normal form \(z_T = a^{-1}b\).  
Therefore, the total number of letters in \(\widetilde{\beta}_m\) is at most \((s+m)n - N_m\).
Since \(\widetilde{\beta}_m\) is positive, its number of letters coincides with its exponent sum \(\xi(\widetilde{\beta}_m)\), which is invariant under conjugation.  
In particular,
\[
\xi(\widetilde{\beta}_m) = \xi(\beta_m) = \xi(z_T) + m n.
\]
Thus,
\[
\xi(z_T) + m n \le (s+m)n - N_m,
\]
which gives
\[
N_m \le s n - \xi(z_T).
\]
The right-hand side is independent of \(m\), as required.

\medskip

To prove that \(\beta_m \in A_Z\), we study the support of \(\widetilde{\beta}_m\) for a suitable value of \(m\), which we now define.  
Recall that in spherical-type Artin groups, \(\Delta_U^2\) is always central in \(A_U\) for any subset \(U\) of standard generators.  
For technical reasons (to be explained below), we will work with this square rather than with \(\Delta_U\) itself.

We have shown that the sequences \(\{c_i\}_{i \geq 1}\) and \(\{R_i\}_{i \geq 1}\) are finite; hence, the subsequences \(\{c_{2i}\}_{i \geq 1}\) and \(\{R_{2i}\}_{i \geq 1}\) are also finite.  
Therefore, there exist integers \(M < m_1 < m_2\) such that
\[
R_{2m_1} = R_{2m_2} =: R, \quad
c_{2m_1} = c_{2m_2} =: c, \quad
U_{2m_1} = U_{2m_2} =: U
\]
(notice that the number of subsets of standard generators is also finite).  
Set \(t = m_2 - m_1\) and \(N = N_{2m_1} = N_{2m_2}\).
To compare \(\widetilde{\beta}_{2m_1}\) and \(\widetilde{\beta}_{2m_2}\), note that
\[
\widetilde{\beta}_{2m_1} = c^{-1} \beta_{2m_1} c,
\]
and
\[
\widetilde{\beta}_{2m_2} 
= c^{-1} \beta_{2m_2} c
= c^{-1} z_T \Delta_Z^{2m_2} c
= c^{-1} \beta_{2m_1} \Delta_Z^{2t} c
= \widetilde{\beta}_{2m_1} \, (c^{-1} \Delta_Z^{2t} c).
\]
On the other hand, since \(\widetilde{\beta}_{2m_2} = \Delta_U^{2m_2-N} R\), we have
\[
\widetilde{\beta}_{2m_2}
= \Delta_U^{2m_2 - 2m_1} \, \Delta_U^{2m_1 - N} R
= \Delta_U^{2t} \, \widetilde{\beta}_{2m_1}.
\]
Here we use the fact that \(\Delta_U^2\) is central in \(A_U\) to deduce
\[
\widetilde{\beta}_{2m_1} \Delta_U^{2t}
= \widetilde{\beta}_{2m_1} \, (c^{-1} \Delta_Z^{2t} c),
\]
which implies
\[
\Delta_U^{2t} = c^{-1} \Delta_Z^{2t} c.
\]
Since the parabolic closure of \(\Delta_Z^{2t}\) is \(A_Z\) and that of \(\Delta_U^{2t}\) is \(A_U\), the solution of \autoref{exo_parabolicclosure} gives
\[
A_U = c^{-1} A_Z c.
\]
Moreover, by the proof of \autoref{theorem:unique_parabolic_closure}, the parabolic closure of the recurrent element \(\widetilde{\beta}_{2m_1}\) is precisely \(A_U\), the standard parabolic subgroup defined by its support.  
Thus,
\[
P_{\beta_{2m_1}} = c \, P_{\widetilde{\beta}_{2m_1}} \, c^{-1} = c A_U c^{-1} = A_Z.
\]
Therefore \(\beta_{2m_1} \in A_Z\), which completes the proof of the theorem.
\end{proof}

\subsubsection{Some infinite cases taking parabolic subgroups as stabilizers in complexes}\label{section:GGT}

The objective of this section is to explain how we can solve the problem of the intersection of parabolic subgroups in some infinite-type cases viewing the parabolic subgroups as stabilizers in cubical or simplicial complexes, using an induction that has as base case the spherical case seen before. 

\paragraph{Spherical-type parabolic subgroups in FC-type Artin groups}

We now explain the approach taken by Rose Morris-Wright in her PhD thesis. She uses the cubical decomposition of the Deligne complex and the CAT(0) property —specifically, the uniqueness of geodesic paths between vertices stabilized by spherical-type parabolic subgroups— to study combinatorial geodesics, that is, minimal-length paths in the 1-skeleton of the complex.

\medskip

Consider the partially ordered set $$ AS^f= \{\alpha A_T, \alpha\in A_S, T\in S^f\} $$ used to define the Deligne complex. It is not difficult to see that every interval \([ \alpha A_{T_1}, \alpha A_{T_2} ]\) spans a cube of dimension \( |T_2 \setminus T_1| \).  For example, take the portion of the Deligne complex corresponding to the Artin group
\[
\langle a,b,c,d \mid aba=bab,\; aca=cac,\; adad=dada,\; bc=cb,\; cdc=dcd \rangle.
\] 
depicted in \autoref{Dibujo:derived}. We can see that $[1\cdot A_{\{a\}},1\cdot A_{\{a,d\}}]$ spans a cube of dimension 1 (an edge), and $[1\cdot A_{\{a\}},1\cdot A_{\{a,b,c\}}]$ spans a cube of dimension 2 (a square with vertices $\{a\},\{a,b\},\{a,c\}$ and $\{a,b,c\}$).
By declaring that each cube of dimension \( m \) is isometric to the cube \( [0,1]^m \), we obtain the cubical decomposition of the Deligne complex. A classical result of \cite{CharneyDavis} shows that, in the FC-type case, this complex is CAT(0).

\medskip

The action of an element \( h \) of an FC-type Artin group on a vertex \( gA_T \) sends it to the vertex \( hgA_T \). This induces an action by isometries on the full complex. The stabilizer of a vertex \( gA_T \) is precisely the spherical-type parabolic subgroup \( gA_T g^{-1} \).

\begin{exercise}\label{exercise:pointwise_geodesic}
If \( g \in A_S\) fixes two vertices of the cubical decomposition of the Deligne complex of an FC-type Artin group $A_S$, then it fixes pointwise any combinatorial geodesic between them. \emph{(See solution in page \pageref{sol:pointwise}.)}
\end{exercise}

We say that an edge in the cubical decomposition of the Deligne complex is \emph{downward} if it goes from a vertex \( gA_X \) to a vertex \( gA_{X \setminus \{x\}} \) for some \( x \in X \), and \emph{upward} if it goes from a vertex \( gA_X \) to a vertex \( gA_{X \cup \{y\}} \) for some \( y \in S \setminus X \). In a combinatorial path (consisting of edges), we say that a vertex is a \emph{turning point} if it lies between an upward and a downward edge.

\begin{theorem}[{\citealp[Theorem 3.1]{Rose}}]\label{theorem:Rose}
In an FC-type Artin group \( A_S \), the intersection of two spherical-type parabolic subgroups \( P \) and \( Q \) is again a parabolic subgroup.
\end{theorem}

\begin{proof}

Pick a vertex \( v_P \) in the cubical decomposition of the Deligne complex stabilized by~\( P \), and a vertex \( v_Q \) stabilized by \( Q \), and let \( p \) be a combinatorial geodesic between them. By \autoref{exercise:pointwise_geodesic}, we know that \( P \cap Q \) fixes~\( p \) pointwise. Since, in a combinatorial path consisting only of upward (or only of downward) edges, the stabilizers of the vertices form a chain of inclusions, we will obtain our result by induction on the number of turning points $v_0=v_P,v_1,\dots,v_m=v_Q$ in \( p \), where we denote by~$P_i$ the (spherical-type) parabolic subgroup that stabilizes the turning point~$v_i$.

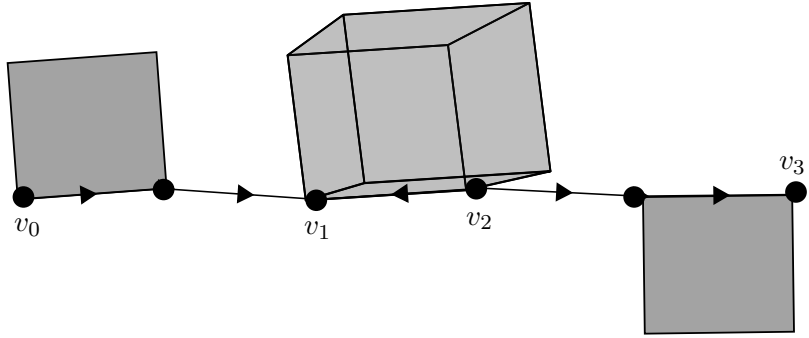
\begin{figure}[h]
  \centering
    \def\svgwidth{0.7\columnwidth} 
\begingroup%
  \makeatletter%
  \providecommand\color[2][]{%
    \errmessage{(Inkscape) Color is used for the text in Inkscape, but the package 'color.sty' is not loaded}%
    \renewcommand\color[2][]{}%
  }%
  \providecommand\transparent[1]{%
    \errmessage{(Inkscape) Transparency is used (non-zero) for the text in Inkscape, but the package 'transparent.sty' is not loaded}%
    \renewcommand\transparent[1]{}%
  }%
  \providecommand\rotatebox[2]{#2}%
  \newcommand*\fsize{\dimexpr\f@size pt\relax}%
  \newcommand*\lineheight[1]{\fontsize{\fsize}{#1\fsize}\selectfont}%
  \ifx\svgwidth\undefined%
    \setlength{\unitlength}{401.19836426bp}%
    \ifx\svgscale\undefined%
      \relax%
    \else%
      \setlength{\unitlength}{\unitlength * \real{\svgscale}}%
    \fi%
  \else%
    \setlength{\unitlength}{\svgwidth}%
  \fi%
  \global\let\svgwidth\undefined%
  \global\let\svgscale\undefined%
  \makeatother%
  \begin{picture}(1,0.39860325)%
    \lineheight{1}%
    \setlength\tabcolsep{0pt}%
    \put(0,0){\includegraphics[width=\unitlength,page=1]{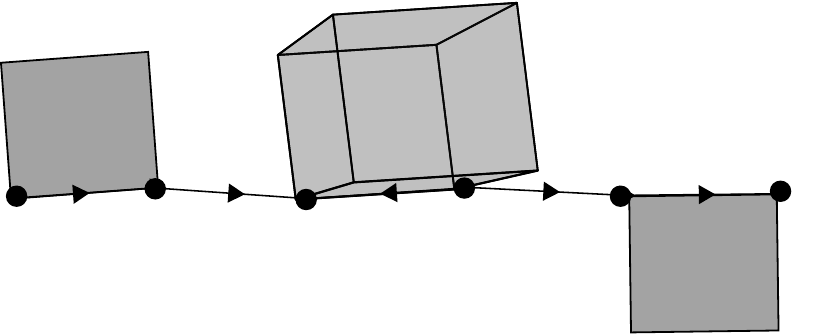}}%
    \put(0.00887979,0.12304764){\color[rgb]{0,0,0}\makebox(0,0)[lt]{\lineheight{1.25}\smash{\begin{tabular}[t]{l}$v_0$\end{tabular}}}}%
    \put(0.35138588,0.11924198){\color[rgb]{0,0,0}\makebox(0,0)[lt]{\lineheight{1.25}\smash{\begin{tabular}[t]{l}$v_1$\end{tabular}}}}%
    \put(0.54271904,0.13064174){\color[rgb]{0,0,0}\makebox(0,0)[lt]{\lineheight{1.25}\smash{\begin{tabular}[t]{l}$v_2$\end{tabular}}}}%
    \put(0.91201425,0.19648937){\color[rgb]{0,0,0}\makebox(0,0)[lt]{\lineheight{1.25}\smash{\begin{tabular}[t]{l}$v_3$\end{tabular}}}}%
  \end{picture}%
\endgroup%

  \medskip
  \caption{An illustration for the proof of \autoref{theorem:Rose}.}\label{Dibujo:pruebaRose}
\end{figure}

If \( m = 0 \), then \( P = Q \), and the result follows. Now suppose that \( P \cap P_{n-1} \) is a spherical-type parabolic subgroup, which by construction fixes pointwise the subpath of \( p \) from \( v_0 = v_P \) to~\( v_{n-1} \). In particular, this implies that \( P \cap Q \subset P \cap P_{n-1} \), and therefore
$
P \cap Q = P \cap Q \cap P_{n-1}.$
We now distinguish two cases:

\begin{enumerate}
\item The final segment of \( p \) from \( v_{n-1} \) to~\( v_n \) consists entirely of upward edges. In this case, \( P_{n-1} \) is strictly contained in \( P_n = Q \), so
\[
P \cap P_{n-1} = P \cap Q \cap P_{n-1} = P \cap Q.
\]

\item The final segment of \( p \) from \( v_{n-1} \) to \( v_n \) consists entirely of downward edges. In this case, \( Q \) is strictly contained in \( P_{n-1} \), and both \( P \cap P_{n-1} \) and \( Q \) are parabolic subgroups contained in the spherical-type parabolic subgroup \( P_{n-1} \). Up to conjugacy, we may assume that \( P_{n-1} \) is standard. Then, by \autoref{theorem_restandardisation} (using just its spherical-type version), both \( P \cap P_{n-1} \) and \( Q \) are parabolic subgroups of \( P_{n-1} \) itself.. Therefore, by \autoref{theorem:intersection}, 
\[
P \cap Q = (P \cap P_{n-1}) \cap Q
\]
is a parabolic subgroup.
\end{enumerate}

\end{proof}

As pointed out in \citep[Remark~3.1]{Rose}, this proof can be adapted to the context of the clique--cube complex. This complex, described implicitly in \citep{GodelleParis} and explicitly in \citep{CharneyMorris}, is constructed in the same way as the Deligne complex, but using the set 
\(\{X \subset S \mid X \text{ is free of } \infty\}\) instead of \(S^f\). 
With this setup, one can prove that if Conjecture~4 holds for every Artin group whose Coxeter graph is complete (under the no-\(\infty\) convention), then in any Artin group the intersection of two parabolic subgroups corresponding to cliques---that is, conjugates of standard parabolic subgroups associated with complete subgraphs---is again a parabolic subgroup.

\paragraph{Some 2-dimensional families including large-type Artin groups}
Here, we describe the approach from \citep{CMV}, which is similar to the previous one but takes place in a simplicial complex known as the Artin complex. In this setting, we study combinatorial geodesics in much the same spirit as in the previous proof. The key property we rely on is \emph{systolicity}, which ensures that if an element acts on the complex fixing two vertices, then it also fixes every combinatorial geodesic connecting them. 
We will also discuss the generalization by \cite{Blufstein} of the result to \((2,2)\)-free 2-dimensional Artin groups—namely, those 2-dimensional Artin groups whose Coxeter graph, in the no-\( \infty \) notation, does not contain two consecutive edges labelled by 2.

\medskip

\begin{figure}[h]
  \centering
    \def\svgwidth{0.4\columnwidth} 
\begingroup%
  \makeatletter%
  \providecommand\color[2][]{%
    \errmessage{(Inkscape) Color is used for the text in Inkscape, but the package 'color.sty' is not loaded}%
    \renewcommand\color[2][]{}%
  }%
  \providecommand\transparent[1]{%
    \errmessage{(Inkscape) Transparency is used (non-zero) for the text in Inkscape, but the package 'transparent.sty' is not loaded}%
    \renewcommand\transparent[1]{}%
  }%
  \providecommand\rotatebox[2]{#2}%
  \newcommand*\fsize{\dimexpr\f@size pt\relax}%
  \newcommand*\lineheight[1]{\fontsize{\fsize}{#1\fsize}\selectfont}%
  \ifx\svgwidth\undefined%
    \setlength{\unitlength}{231.32176677bp}%
    \ifx\svgscale\undefined%
      \relax%
    \else%
      \setlength{\unitlength}{\unitlength * \real{\svgscale}}%
    \fi%
  \else%
    \setlength{\unitlength}{\svgwidth}%
  \fi%
  \global\let\svgwidth\undefined%
  \global\let\svgscale\undefined%
  \makeatother%
  \begin{picture}(1,0.78739844)%
    \lineheight{1}%
    \setlength\tabcolsep{0pt}%
    \put(0,0){\includegraphics[width=\unitlength,page=1]{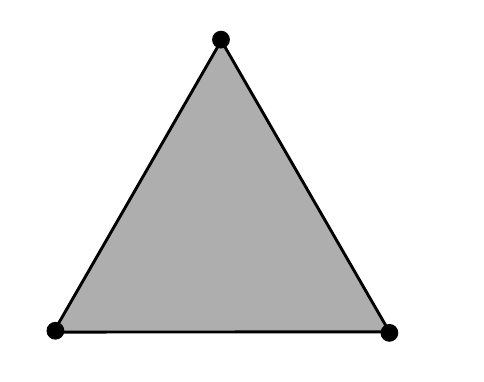}}%
    \put(0.42053798,0.74081416){\color[rgb]{0,0,0}\makebox(0,0)[lt]{\lineheight{1.25}\smash{\begin{tabular}[t]{l}$\{a\}$\end{tabular}}}}%
    \put(0.00125218,0.08276338){\color[rgb]{0,0,0}\makebox(0,0)[lt]{\lineheight{1.25}\smash{\begin{tabular}[t]{l}$\{b\}$\end{tabular}}}}%
    \put(0.8397561,0.08310253){\color[rgb]{0,0,0}\makebox(0,0)[lt]{\lineheight{1.25}\smash{\begin{tabular}[t]{l}$\{c\}$\end{tabular}}}}%
    \put(0.14034135,0.41999838){\color[rgb]{0,0,0}\makebox(0,0)[lt]{\lineheight{1.25}\smash{\begin{tabular}[t]{l}$\{a,b\}$\end{tabular}}}}%
    \put(0.63742311,0.42239076){\color[rgb]{0,0,0}\makebox(0,0)[lt]{\lineheight{1.25}\smash{\begin{tabular}[t]{l}$\{a,c\}$\end{tabular}}}}%
    \put(0.40098004,0.01447858){\color[rgb]{0,0,0}\makebox(0,0)[lt]{\lineheight{1.25}\smash{\begin{tabular}[t]{l}$\{b,c\}$\end{tabular}}}}%
    \put(0.36405286,0.32149849){\color[rgb]{0,0,0}\makebox(0,0)[lt]{\lineheight{1.25}\smash{\begin{tabular}[t]{l}$\{a,b,c\}$\end{tabular}}}}%
  \end{picture}%
\endgroup%

  \medskip
  \caption{The portion of an Artin complex  of an Artin group on three generators $a,b$ and $c$ corresponding to the cosets of the form $1\cdot A_X$, $X\subset S$.}\label{Dibujo:ArtinComplex}
\end{figure}
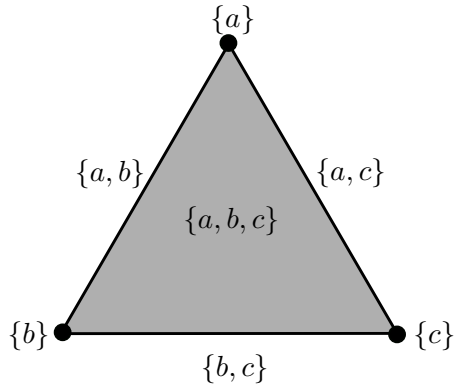

We consider the Artin complex of a large-type Artin group \( A_S \) with \(|S| = n\). 
Recall that the Artin complex is the complex whose barycentric subdivision is the derived complex of $A\mathcal{S}$. Under this viewpoint, cosets associated 
to subsets of \(S\) of size \(k\) correspond to simplices of dimension \(n-1\). See \autoref{Dibujo:ArtinComplex} for an illustration. 

As before, the action of the group on the complex comes from the action of cosets by left multiplication. This action is without inversion. In this complex, the stabilizer of a \( k \)-dimensional simplex is a parabolic subgroup of rank \( n - k -1\). Moreover, for every parabolic subgroup, there exists a (non-unique) simplex whose stabilizer is that subgroup~\cite[Theorem~11]{CMV}. In our example of \autoref{Dibujo:ArtinComplex}, using the sets of generators as labels of the simplices, we have
\[
\begin{array}{lll}
\mathrm{Stab}(\{a\}) = A_{\{b,c\}}, &
\mathrm{Stab}(\{b\}) = A_{\{a,c\}}, &
\mathrm{Stab}(\{c\}) = A_{\{a,b\}}, \\[6pt]

\mathrm{Stab}(\{a,b\}) = A_{\{c\}}, &
\mathrm{Stab}(\{b,c\}) = A_{\{a\}}, &
\mathrm{Stab}(\{a,c\}) = A_{\{b\}}, \\[6pt]

\multicolumn{3}{c}{
\mathrm{Stab}(\{a,b,c\}) = 1.
}
\end{array}
\]
\medskip
Let \(K\) be a simplicial complex and let \(D \subseteq K\) be a simplex. 
The \emph{link} of \(D\) in \(K\), denoted~\(\mathrm{Lk}_K(D)\), 
is the subcomplex of \(K\) that consists of all simplices of \(K\) that are disjoint from~\(D\) but together with~\(D\) span a simplex of~\(K\). We say that a cycle $p$ in the 1-skeleton of the complex is \emph{full} if the only sets of vertices that span a simplex are consecutive vertices on the cycle. We denote the number of edges in $p$ by $|p|$.

\begin{definition}
The \textit{systole} of a simplicial complex \( K \) is defined as
\[
\operatorname{sys}(K) \coloneqq \min \{ |p| \mid p \text{ is an embedded full cycle in } K \} \in \{3, 4, \ldots, \infty\}.
\]
For \( k \in \{3, \ldots, \infty\} \), we say that \( K \) is \textit{locally \( k \)-large} if \( \operatorname{sys}(\operatorname{Lk}_K(D)) \geq k \) for all simplices \( D \subseteq K \). We say that \( K \) is \textit{\( k \)-large} if it is locally \( k \)-large and \( \operatorname{sys}(K) \geq k \). The complex \( K \) is said to be \textit{\( k \)-systolic} if it is connected, simply connected, and locally \( k \)-large. Finally, \( K \) is called \textit{systolic} if it is 6-systolic.
\end{definition}

It was shown in~\cite[Theorem~8]{CMV} that for large-type Artin groups of rank greater than~2, the Artin complex is systolic. As a consequence, any element of a group acting without inversion on the complex that fixes two vertices must also fix every combinatorial geodesic between them~\cite[Lemma~14]{CMV}. We call this property the \emph{two-vertex rigidity} property. Later, by assigning suitable weights to the edges of the complex, \cite{Blufstein} established that for 2-dimensional $(2,2)$-free Artin groups the Artin complex satisfies a more general property, which he called \emph{systolic-by-function}. In this setting, the same two-vertex rigidity property also holds \citep[Theorem~1.2]{Blufstein}.

\begin{exercise}\label{exercise:simplices}
Let $A$ be an Artin group and $K$ be its Artin complex. Prove that if \( P \) and \( Q \) are parabolic subgroups, corresponding to the stabilizers of simplices \( D_1 \) and \( D_2 \) of $K$ that share a face \( F \), then, up to conjugacy, \( P \) and \( Q \) are parabolic subgroups of \( \mathrm{Stab}(F) \). \emph{(See solution in page~\pageref{sol:simplices}.)}
\end{exercise}

 We now turn to the intersection result. We will denote a combinatorial path in the Artin complex by its consecutive edges \( e_1, \dots, e_k \). We say that a standard parabolic subgroup $A_X$ \emph{satisfies the intersection property} if for any two $P$ and $Q$ parabolic subgroups of $A_X$ we have that $P\cap Q$ is a parabolic subgroup.

\begin{lemma}[{Claim 1 in the proof of \citealp[Theorem~11]{CMV}}]\label{lemma:path}
Consider an Artin group \( A_S \)  and its Artin complex. Let \( e_1, \dots, e_k \) be a combinatorial path in the complex. If the proper standard parabolic subgroups of \( A_S \) satisfy the intersection property, then the intersection of the edge stabilizers is the stabilizer of a simplex that contains~\( e_k \). That is,
\[
\bigcap_{1 \leq i \leq k} \mathrm{Stab}(e_i) = \mathrm{Stab}(D), \quad e_k \subseteq D.
\]
\end{lemma}

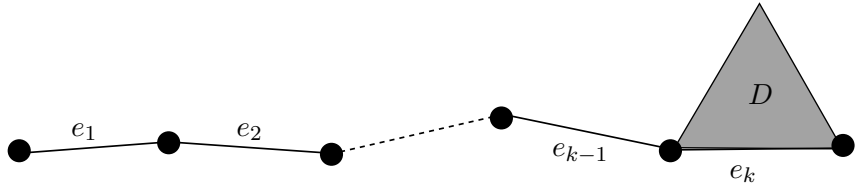
\begin{figure}[h]
  \centering
    \def\svgwidth{0.7\columnwidth} 
\begingroup%
  \makeatletter%
  \providecommand\color[2][]{%
    \errmessage{(Inkscape) Color is used for the text in Inkscape, but the package 'color.sty' is not loaded}%
    \renewcommand\color[2][]{}%
  }%
  \providecommand\transparent[1]{%
    \errmessage{(Inkscape) Transparency is used (non-zero) for the text in Inkscape, but the package 'transparent.sty' is not loaded}%
    \renewcommand\transparent[1]{}%
  }%
  \providecommand\rotatebox[2]{#2}%
  \newcommand*\fsize{\dimexpr\f@size pt\relax}%
  \newcommand*\lineheight[1]{\fontsize{\fsize}{#1\fsize}\selectfont}%
  \ifx\svgwidth\undefined%
    \setlength{\unitlength}{376.22698975bp}%
    \ifx\svgscale\undefined%
      \relax%
    \else%
      \setlength{\unitlength}{\unitlength * \real{\svgscale}}%
    \fi%
  \else%
    \setlength{\unitlength}{\svgwidth}%
  \fi%
  \global\let\svgwidth\undefined%
  \global\let\svgscale\undefined%
  \makeatother%
  \begin{picture}(1,0.2167432)%
    \lineheight{1}%
    \setlength\tabcolsep{0pt}%
    \put(0,0){\includegraphics[width=\unitlength,page=1]{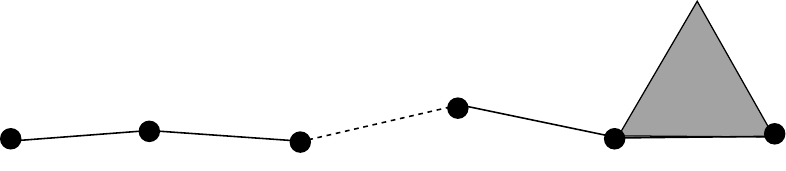}}%
    \put(0.87420004,0.09630275){\color[rgb]{0,0,0}\makebox(0,0)[lt]{\lineheight{1.25}\smash{\begin{tabular}[t]{l}$D$\end{tabular}}}}%
    \put(0.07514573,0.05853199){\color[rgb]{0,0,0}\makebox(0,0)[lt]{\lineheight{1.25}\smash{\begin{tabular}[t]{l}$e_1$\end{tabular}}}}%
    \put(0.27023783,0.0584594){\color[rgb]{0,0,0}\makebox(0,0)[lt]{\lineheight{1.25}\smash{\begin{tabular}[t]{l}$e_2$\end{tabular}}}}%
    \put(0.64305095,0.03450132){\color[rgb]{0,0,0}\makebox(0,0)[lt]{\lineheight{1.25}\smash{\begin{tabular}[t]{l}$e_{k-1}$\end{tabular}}}}%
    \put(0.85223528,0.00890213){\color[rgb]{0,0,0}\makebox(0,0)[lt]{\lineheight{1.25}\smash{\begin{tabular}[t]{l}$e_{k}$\end{tabular}}}}%
  \end{picture}%
\endgroup%

  \medskip
  \caption{An illustration for \autoref{lemma:path}.}\label{Dibujo:pruebalarge}
\end{figure}

\begin{proof}
We proceed by induction on \( k \). If \( k = 1 \), the proof is trivial. Now suppose that the result holds for \( k - 1 \). Then we have
\[
\bigcap_{1 \leq i \leq k} \mathrm{Stab}(e_i) = \mathrm{Stab}(D') \cap \mathrm{Stab}(e_{k+1}), \quad e_{k-1} \subseteq D'.
\]
Let \( v \) be a vertex contained in both \( e_{k-1} \) and \( e_k \). Up to conjugacy, we may assume that \( \mathrm{Stab}(v) \) is standard, hence isomorphic to an Artin group on \( n-1 \) generators that satisfies the intersection property. By \autoref{exercise:simplices}, this implies that both \( \mathrm{Stab}(D') \) and \( \mathrm{Stab}(e_k) \) are parabolic subgroups of \( \mathrm{Stab}(v) \). Thus, \( \mathrm{Stab}(D') \cap \mathrm{Stab}(e_k) \) is a parabolic subgroup of \( \mathrm{Stab}(v) \) contained in \( \mathrm{Stab}(e_k) \), so it is a parabolic subgroup of \( \mathrm{Stab}(e_k) \) by \autoref{theorem_restandardisation}. Geometrically, \( \mathrm{Stab}(D') \cap \mathrm{Stab}(e_k) \) is the stabilizer of some simplex containing~\( e_k \).
\end{proof}

\begin{theorem}
In an Artin group \( A_S \)  with an Artin complex satisfying the two-vertex rigidity property for $|S|>2$, the intersection of two parabolic subgroups \( P \) and \( Q \) is again a parabolic subgroup.
\end{theorem}

\begin{proof}
For the case \( |S| = 2 \), if the relation between the two generators exists—that is, it is different from \( \infty \)—then the group is of spherical type, and the result follows from \autoref{theorem:intersection}. Otherwise, it is an FC-type Artin group with only proper spherical-type parabolic subgroups, and the result follows from \autoref{theorem:Rose}.

For the case \( |S| = n > 2 \), consider the action of the group on its Artin complex. We proceed by induction on $n$, assuming that the result holds for \( n-1 \), so the hypotheses of \autoref{lemma:path} apply. Let \( D_1 \) and \( D_2 \) be simplices whose stabilizers are \( P \) and \( Q \), respectively. We define a combinatorial path $p=e_1,\dots, e_k$ that travels along very vertex of $D_1$ then follows a combinatorial geodesic between this last vertex and a vertex of $D_2$ and then travels along every vertex of $D_2$. Since the action of the group on the complex is without inversion,  all the vertices in $D_1$ and $D_2$ are fixed by $P\cap Q$, so it fixes the portion of the path that is a combinatorial geodesic between a vertex of $D_1$ and a vertex of $D_2$, hence fixing the whole $p$. Also, the action without inversion implies that the stabilizer of a simplex is the set the fixes its set of vertices. In particular, $P\cap Q$ will be the set of elements fixing the vertices of both $D_1$ and $D_2$. Combining this with \autoref{lemma:path} we obtain $$P\cap Q = \mathrm{Stab}(D_1)\cap \mathrm{Stab}(D_2) = \bigcap_{1 \leq i \leq k} \mathrm{Stab}(e_i)=  \mathrm{Stab}(D), \quad e_k \subseteq D.$$
\end{proof}

\begin{corollary}
The intersection of two parabolic subgroup in a (2,2)-free 2 dimensional Artin group (which includes large Artin groups) is again a parabolic subgroup. 
\end{corollary}

\subsubsection{Some FC cases using Bass-Serre theory}\label{section:Bass-Serre}

\cite{Moller} extended \autoref{theorem:Rose} by employing, instead of the Deligne complex or the clique complex, a clever use of Bass–Serre theory. Here, we present the proof in the FC-type case, even though the original result is more general, as we will briefly mention the general case at the end of this section.

\medskip

Let \(G_1\) and \(G_2\) be groups and let \(H\) be a group with injective homomorphisms 
\[
\iota_1 : H \hookrightarrow G_1, 
\qquad 
\iota_2 : H \hookrightarrow G_2.
\]
 If
\[
G_1 = \langle S_1 \mid R_1 \rangle,
\qquad 
G_2 = \langle S_2 \mid R_2 \rangle,
\qquad 
H = \langle T \mid R_H \rangle,
\]
and the embeddings \(\iota_1, \iota_2\) send each generator \(t \in T\) to words 
\(w_t^{(1)} \in S_1^{\ast}\) and \(w_t^{(2)} \in S_2^{\ast}\), then the \emph{amalgamated product} of \(G_1\) and \(G_2\) over \(H\) is the group
\[
G_1 \ast_H G_2
= \big\langle S_1 \sqcup S_2 
\;\big|\; R_1,\, R_2,\, w_t^{(1)} = w_t^{(2)} \ \text{for all } t \in T \big\rangle.
\]
That is, one takes the free product \(G_1 \ast G_2\) and imposes the relations identifying the two images of \(H\) inside \(G_1\) and \(G_2\).

\begin{exercise}\label{exo:FC-decomposition}
Let \( A_\Gamma = A_S \) be an Artin group with associated Coxeter graph \( \Gamma \) (with the no-$\infty$ convention), and let~\( \Gamma_1 \) and~\( \Gamma_2 \) be two subgraphs of \( \Gamma \) such that \( \Gamma = \Gamma_1 \cup \Gamma_2 \) and \( \Gamma_1 \cap \Gamma_2 \neq \emptyset \). Show that \( A_\Gamma \) admits the following amalgamated product decomposition:
\[
A_\Gamma = A_{\Gamma_1} *_{A_{\Gamma_1 \cap \Gamma_2}} A_{\Gamma_2}.
\]
In particular, if \( m_{s,t} = \infty \) for some \( s,t \in S \), then:
\[
A_S = A_{S \setminus \{s\}} *_{A_{S \setminus \{s,t\}}} A_{S \setminus \{t\}}.
\]
This shows that if \( A_S \) is of FC type, then it can be expressed as an amalgamated product of spherical-type Artin groups. \emph{(See solution in page~\pageref{sol:FC}.)}
\end{exercise}

Given an amalgamated product of groups \( G = G_1 *_H G_2 \), we construct its Bass–Serre tree as follows. There are two types of vertices: those corresponding to the cosets \( gG_1 \) and those corresponding to the cosets \( gG_2 \), for \( g \in G \). The edges correspond to the cosets \( gH \), and an edge \( cH \) connects two vertices \( aG_1 \) and \( bG_2 \) precisely when \( aG_1 \cap bG_2 = cH \). The group \( G \) acts simplicially on this tree by left multiplication, and the stabilizers of the vertices and edges are, respectively, \( aG_1a^{-1} \), \( bG_2b^{-1} \), and \( cHc^{-1} \).

The construction of the Artin complex discussed earlier—and more generally, the theory of complexes of groups—can be viewed as a generalization of Bass–Serre theory. For more details about Bass-Serre theory, we refer the reader to Serre’s classic book~\citep{Serre1,Serre2}.

\begin{theorem}\label{theorem:intersection:Bass-Serre}
Let \( A_S \) be an FC-type Artin group. Then the intersection of a spherical-type parabolic $P$ subgroup with any other parabolic subgroup $Q$ is again a parabolic subgroup.
\end{theorem}

\begin{proof}
We know that we can write \( P = gA_Xg^{-1} \) and \( Q = hA_Yh^{-1} \), where \( A_X \) is of spherical type. We prove the result by induction on the number \( k \) of pairs \( s,t \in S \) such that \( m_{s,t} = \infty \). 

If \( k = 0 \), then both \( P \) and \( Q \) are of spherical type, and the result follows from \autoref{theorem:Rose}.
Now suppose the result holds for \( k-1 \), and choose a pair \( s,t \in S \) with \( m_{s,t} = \infty \). By \autoref{exo:FC-decomposition}, we know that
\[
A_S = A_{S \setminus \{s\}} *_{A_{S \setminus \{s,t\}}} A_{S \setminus \{t\}}.
\]
Set \( I = S \setminus \{s\} \), \( J = S \setminus \{t\} \), and \( K = S \setminus \{s,t\} \). Let \( T \) be the Bass–Serre tree associated to this amalgam. Since \( A_X \) is of spherical type, it contains no pair with \( m_{s,t} = \infty \), and hence \( X \subset I \) or \( X \subset J \). Thus, there exists a vertex \( u := a_0U_0 \) of \( T \), with \( U_0 \in \{I, J\} \), such that \( P = gA_Xg^{-1} \subset a_0A_{U_0}a_0^{-1} = \mathrm{Stab}(u) \). From here, we distinguish two cases:

\begin{enumerate}
\item \( \{s,t\} \not\subset Y \). In this case, \( Y \subset I \) or \( Y \subset J \), so there is a vertex \( v := b_0V_0 \) in \( T \) such that \( Q \subset b_0A_{V_0}b_0^{-1} = \mathrm{Stab}(v) \). 

We proceed by induction on the distance \( d \) between \( u \) and \( v \). Let \( p \) be the unique geodesic between them, with consecutive vertices:
\[
u_0 = u, \ u_1 = a_1U_1, \ \dots, \ u_d = v,
\]
and denote by \( e_i = c_iA_K \) the edge connecting \( u_{i-1} \) and \( u_i \). 

If \( d = 0 \), then both \( P \) and \( Q \) are contained in \( a_0A_{U_0}a_0^{-1} \). Up to conjugation, we may assume \( a_0 = 1 \), so the group lies in an Artin group with at most \( k - 1 \) $\infty$'s. We apply the inductive hypothesis on \( k \).

Now suppose the result holds for distance \( d - 1 \). Since \( P \subset \mathrm{Stab}(u) \) and \( Q \subset \mathrm{Stab}(v) \), the intersection \( P \cap Q \) fixes the path \( p \), so in particular it fixes all vertices \( u_i \) and all edges~\( e_i \). 
Up to conjugation by \( a_0 \), we may assume \( P \) and \( c_1A_Kc_1^{-1} \) are parabolic subgroups of~\( A_{U_0} \). By the inductive hypothesis on~\( k \), their intersection \( P_1 := P \cap c_1A_Kc_1^{-1} \) is a parabolic subgroup. Since it is contained in \( P \), it is of spherical type, and so are all of its conjugates, so we may conjugate back by \( a_0^{-1} \). We then have:
\[
P \cap Q = P \cap (c_1A_Kc_1^{-1}) \cap Q = P_1 \cap Q.
\]
Now, since \( P_1 \subset c_1A_Kc_1^{-1} \), it stabilizes the edge \( e_1 \), and therefore also stabilizes \( u_1 \). So \( P_1 \subset a_1A_{U_1}a_1^{-1} \), and the problem reduces to a path of length \( d - 1 \). We may now apply the inductive hypothesis on \( d \).

\item \( \{s,t\} \subset Y \). By \autoref{exo:FC-decomposition}, we have:
\[
A_Y = A_{Y \setminus \{s\}} *_{A_{Y \setminus \{s,t\}}} A_{Y \setminus \{t\}}.
\]
Set \( Y_I = Y \setminus \{s\} \), \( Y_J = Y \setminus \{t\} \), and \( Y_K = Y \setminus \{s,t\} \), and let \( T_Y \) be the Bass–Serre tree of this decomposition. 

We claim that there is an embedding of \( T_Y \) into \( T \) such that every vertex \( gA_{Y_U} \in T_Y \) maps to \( gA_U \in T \), where \( U \in \{I, J, K\} \). To see this, we need to show that for all \( g_1, g_2 \in A_Y \),
\[
g_1A_{Y_U} = g_2A_{Y_U} \quad \Leftrightarrow \quad g_1A_U = g_2A_U.
\]
Using the result of \cite{Vanderlek} on intersections of standard parabolic subgroups, we know \( A_Y \cap A_U = A_{Y \cap U} = A_{Y_U} \), so:
\[
g_1A_{Y_U} = g_2A_{Y_U} \quad \Leftrightarrow \quad g_1^{-1}g_2 \in A_{Y_U} = A_Y \cap A_U.
\]
Since $g_1^{-1}g_2$ belongs to $A_Y$, the later is equivalent to
$$g_1^{-1}g_2\in A_U\Leftrightarrow g_1A_{U}=g_2A_{U}. $$

As \( Q = hA_Yh^{-1} \), we consider the translated subtree \( hT_Y \subset T \), and examine the path of minimal length from \( u \) to a vertex of \( hT_Y \). Let \( v = hgA_V \), with \( V \in \{I, J\} \) and \( g \in A_Y \), be the unique vertex of \( hT_Y \) at minimal distance from \( u \) (see \autoref{Dibujo:arbol}). Since \( P \subset \mathrm{Stab}(u) \) and \( Q \subset \mathrm{Stab}(hT_Y) \), we know \( P \cap Q \subset \mathrm{Stab}(v) \), so:
\[
P \cap Q = P \cap (hgA_Vg^{-1}h^{-1}) \cap Q.
\]
Note that \( Q = hA_Yh^{-1} = hgA_Yg^{-1}h^{-1} \) for any \( g \in A_Y \), so up to conjugation by \( hg \), this reduces to the intersection \( A_V \cap A_Y = A_{V \cap Y} \). Therefore:
\[
P \cap Q = P \cap (hgA_{V \cap Y}g^{-1}h^{-1}).
\]
Since \( \{s,t\} \not\subset V \cap Y \), we may now apply Case 1 to conclude the result.

\medskip

\begin{figure}[H]
  \centering
  \def\svgwidth{0.8\columnwidth} 
\begingroup%
  \makeatletter%
  \providecommand\color[2][]{%
    \errmessage{(Inkscape) Color is used for the text in Inkscape, but the package 'color.sty' is not loaded}%
    \renewcommand\color[2][]{}%
  }%
  \providecommand\transparent[1]{%
    \errmessage{(Inkscape) Transparency is used (non-zero) for the text in Inkscape, but the package 'transparent.sty' is not loaded}%
    \renewcommand\transparent[1]{}%
  }%
  \providecommand\rotatebox[2]{#2}%
  \newcommand*\fsize{\dimexpr\f@size pt\relax}%
  \newcommand*\lineheight[1]{\fontsize{\fsize}{#1\fsize}\selectfont}%
  \ifx\svgwidth\undefined%
    \setlength{\unitlength}{246.03429839bp}%
    \ifx\svgscale\undefined%
      \relax%
    \else%
      \setlength{\unitlength}{\unitlength * \real{\svgscale}}%
    \fi%
  \else%
    \setlength{\unitlength}{\svgwidth}%
  \fi%
  \global\let\svgwidth\undefined%
  \global\let\svgscale\undefined%
  \makeatother%
  \begin{picture}(1,0.53530969)%
    \lineheight{1}%
    \setlength\tabcolsep{0pt}%
    \put(0,0){\includegraphics[width=\unitlength,page=1]{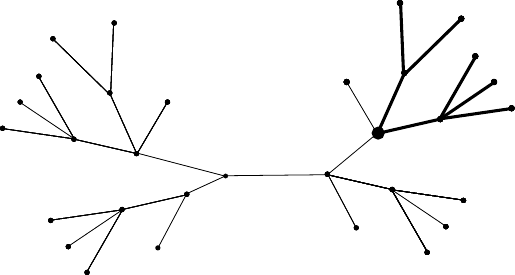}}%
    \put(0.01526898,0.35344131){\color[rgb]{0,0,0}\makebox(0,0)[lt]{\lineheight{1.25}\smash{\begin{tabular}[t]{l}$u$\end{tabular}}}}%
    \put(0.77698901,0.31396534){\color[rgb]{0,0,0}\makebox(0,0)[lt]{\lineheight{1.25}\smash{\begin{tabular}[t]{l}$hT_Y$\end{tabular}}}}%
    \put(0,0){\includegraphics[width=\unitlength,page=2]{arbol.pdf}}%
    \put(0.74293084,0.24275596){\color[rgb]{0,0,0}\makebox(0,0)[lt]{\lineheight{1.25}\smash{\begin{tabular}[t]{l}$v$\end{tabular}}}}%
  \end{picture}%
\endgroup%

  \caption{An illustration for Case 2 on the proof of \autoref{theorem:intersection:Bass-Serre}.}\label{Dibujo:arbol}
\end{figure}

\end{enumerate}
\end{proof}

Note that one can also use \autoref{exo:FC-decomposition} and apply the exact same proof to show that if the intersection property for parabolic subgroups holds in Artin groups whose Coxeter graphs are complete, then the intersection of a parabolic subgroup corresponding to a complete subgraph with any other parabolic subgroup is again a parabolic subgroup. This extends the result of Morris-Wright mentioned after \autoref{theorem:Rose}.

\subsubsection{Artin groups admitting algebraic retractions}\label{section:retractions}

A more recent strategy for studying the structure of parabolic subgroups—particularly in cases where previous results do not apply—is based on the use of algebraic retractions. Given a group homomorphism \( \rho \colon G \to H \), where \( H \leq G \), we say that \( \rho \) is a \emph{retraction} if the restriction \( \rho|_H \) is the identity map on \( H \).

\medskip

In the setting of Artin groups with only even labels, such retractions arise naturally. Indeed, for any subset \( X \subseteq S \), there exists a canonical retraction \( \rho_X \colon A_S \to A_X \) defined by
\[
\rho_X(s) = 
\begin{cases}
s & \text{if } s \in X, \\
1 & \text{otherwise}.
\end{cases}
\]
This map extends to a well-defined group homomorphism due to the parity of the defining relations.

\medskip

\cite{AntolinFoniqi} exploit these natural retractions to analyse parabolic subgroups in even FC-type Artin groups. Among other results, they show that the intersection of two parabolic subgroups is again parabolic in this setting.

\medskip

As an illustration of how retractions simplify arguments, consider the following proposition, originally stated for general Artin groups in~\citep*[Proposition~2.6]{Moller}. The proof below, adapted from~\citep[Lemma~3.4]{AntolinFoniqi}, is significantly more direct when retractions are available.

\begin{proposition}\label{lem: containing implies equality}
Let \( A_S \) be an Artin group that admits retractions, and let \( g, h \in A_S \) and \( X \subseteq S \). If
\[
gA_Xg^{-1} \leqslant hA_Xh^{-1},
\]
then \( gA_Xg^{-1} = hA_Xh^{-1} \).
\end{proposition}

\begin{proof}
The inclusion \( gA_Xg^{-1} \leqslant hA_Xh^{-1} \) is equivalent to
\[
h^{-1}gA_Xg^{-1}h \leqslant A_X.
\]
Applying the retraction \( \rho_X \) to both sides yields
\[
\rho_X(h^{-1}gA_Xg^{-1}h) = h^{-1}gA_Xg^{-1}h,
\]
since \( \rho_X \) is the identity on any subgroup of \( A_X \). But by construction, \( \rho_X(h^{-1}g) \in A_X \), so the conjugate collapses to \( A_X \). Therefore,
\[
h^{-1}gA_Xg^{-1}h = A_X,
\]
and conjugating back, we obtain \( gA_Xg^{-1} = hA_Xh^{-1} \), as claimed.
\end{proof}

In what follows, we will not reproduce the full argument of~\cite{AntolinFoniqi} for FC-type even Artin groups. Instead, we will show that, with minimal additional effort, retractions can be used to compute the intersection of certain pairs of parabolic subgroups in even Artin groups that are not of FC-type that admit retractions—cases that were previously out of reach using other techniques.

\begin{lemma}[{\citealp[Lemma~3.3]{AntolinFoniqi}}]\label{lem:intersection_of_two_parabolic_subgroups_simple}
Let \( A_S \) be an even Artin group. For any \( f, g \in A_S \) and \( X, Y \subseteq S \), there exist \( x \in A_X \) and \( y \in A_Y \) such that
\[
fA_Xf^{-1} \cap gA_Yg^{-1} = fxA_{X \cap Y}x^{-1}f^{-1} \cap gyA_{X \cap Y}y^{-1}g^{-1}.
\]
\end{lemma}

\begin{proof}
We start with the identity:
\[
fA_Xf^{-1} \cap gA_Yg^{-1} = f[A_X \cap (f^{-1}g)A_Y(f^{-1}g)^{-1}]f^{-1}.
\]
Set \( h = f^{-1}g \), and consider \( P = A_X \cap hA_Yh^{-1} \). By the assumption that retractions exist and using the standard intersection result \( A_X \cap A_Y = A_{X \cap Y} \) of Van der Lek, we compute:
\begin{align*}
P = \rho_X(P) &= \rho_X(A_X \cap hA_Yh^{-1}) \\
&\subseteq \rho_X(A_X) \cap \rho_X(hA_Yh^{-1}) \\
&= A_X \cap \rho_X(h)\rho_X(A_Y)\rho_X(h^{-1}) \\
&= \rho_X(h)A_{X \cap Y}\rho_X(h)^{-1}.
\end{align*}
Set \( x = \rho_X(h) \in A_X \). Then \( P \subseteq xA_{X \cap Y}x^{-1} \subseteq A_X \). But since \( P = A_X \cap hA_Yh^{-1} \), we also have:
\[
P = hA_Yh^{-1} \cap xA_{X \cap Y}x^{-1}.
\]
Now conjugate this equation by \( h^{-1} \), and define \( P' = h^{-1}Ph \) and \( k = h^{-1}x \). Then:
\[
P' = A_Y \cap kA_{X \cap Y}k^{-1}.
\]
Applying the same reasoning to \( P' \) using the retraction \( \rho_Y \), we obtain:
\[
P' \subseteq \rho_Y(k)A_{X \cap Y}\rho_Y(k)^{-1}.
\]
Set \( y = \rho_Y(k) \in A_Y \), so that \( P' \subseteq yA_{X \cap Y}y^{-1} \subseteq A_Y \). Combining this with the expression for~\( P' \), we conclude:
\[
P = xA_{X \cap Y}x^{-1} \cap hyA_{X \cap Y}y^{-1}h^{-1}.
\]
Returning to the original setup, we get:
\[
fA_Xf^{-1} \cap gA_Yg^{-1} = fPf^{-1} = fxA_{X \cap Y}x^{-1}f^{-1} \cap gyA_{X \cap Y}y^{-1}g^{-1},
\]
as claimed.
\end{proof}


%

\medskip

Now, let us take the even Artin group given by the following Coxeter graph with the no-2 convention: 

\begin{center}
\begin{tikzpicture}\label{grafotikz}

\tikzstyle{vertex}=[circle, draw, fill=black, inner sep=0pt, minimum size=6pt]

\node[vertex] (A) at (0,0) {};
\node[vertex] (B) at (2,0) {};
\node[vertex] (C) at (4,0) {};
\node[vertex] (D) at (6,0) {};
\node[vertex] (E) at (4,1.5) {};
\node[vertex] (F) at (6,1.5) {};

\draw (A) -- node[below] {4} (B);
\draw (C) -- node[below] {6} (D);
\draw (C) -- node[left] {$\infty$} (E);
\draw (D) -- node[right] {8} (F);
\draw (E) -- node[above] {$\infty$} (F);

\node[above=0.1cm] at (A) {a};
\node[above=0.1cm] at (B) {b};
\node[below=0.1cm] at (C) {c};
\node[below=0.1cm] at (D) {d};
\node[left=0.1cm] at (E) {e};
\node[right=0.1cm] at (F) {f};

\end{tikzpicture}
\end{center}

\medskip

This Artin group is not of FC type. However, given two parabolic subgroups \( P = gA_Xg^{-1} \) and \( Q = hA_Yh^{-1} \), we know from \autoref{lem:intersection_of_two_parabolic_subgroups_simple} that if \( X \cap Y = \emptyset \)—for instance, \( X = \{a, f\} \) and \( Y = \{b, c, e\} \)—then \( P \cap Q = \{1\} \).

Moreover, the lemma tells us that if \( X \cap Y \) lies entirely within a direct component where the intersection of parabolic subgroups is already understood, then the problem reduces to computing the intersection within that component. For example, if \( X = \{a, c\} \) and \( Y = \{a, d\} \), then \( P \cap Q \) is the intersection of two cyclic parabolic subgroups inside the dihedral subgroup~\( A_{\{a, b\}} \), which can be computed using \autoref{theorem:intersection}.

Following the proof of the lemma further, one can observe that if \( X \) contains an entire direct component~\( X' \) and \( X \cap Y \subseteq X' \), then \( P \cap Q \) is a parabolic subgroup supported on a subset of~\( X' \). We leave the following specific instance as an exercise so the reader can develop a proof:

\begin{exercise}\label{exo:retracto}
Consider the Artin group \( A_S \) described above. Prove that if \( X = \{a, c, d, e, f\} \) and \( Y = \{b, c, d\} \), then for any \( g, h \in A_S \), the intersection \( gA_Xg^{-1} \cap hA_Yh^{-1} \) is a parabolic subgroup over \( \{c, d\} \), contained in \( A_{\{c, d, e, f\}} \). \emph{(See solution in page~\pageref{sol:retracto}.)}
\end{exercise}

Beyond the even case, in~\citep{CisnerosCumplidoFoniqi} we provided a complete classification of FC-type Artin groups that admit retractions. The same authors, jointly with Luis Paris, extend this to a classification of all Artin groups admitting retractions \citep{CCFL}. The techniques developed by Antolín and Foniqi are highly specific to even FC-type Artin groups. Although in~\citep{CisnerosCumplidoFoniqi} we managed to extend some of their results to the broader class of FC-type Artin groups that admit retractions—and were able to compute the intersection of certain pairs of parabolic subgroups, as discussed above—we were not able to prove that the intersection of every pair of parabolic subgroups is again parabolic. This question therefore remains open.

\subsubsection{Euclidean braids}\label{section:euclidean} 

Another approach to studying certain Artin groups involves constructing morphisms to other groups that contain well-understood Artin subgroups. These target groups serve as a framework where known structural results can be transferred back to the original group via the morphism. This strategy was used, for example, by \cite{CalvezCisneros} to show that the complex of irreducible parabolic subgroups is Gromov-hyperbolic for the Artin group of spherical type $B_n$ and for the Euclidean braid group. To illustrate the effectiveness of this method, we explain the proof in \citep{CGP}, to show that the intersection of two parabolic subgroups in the Euclidean braid group is again a parabolic subgroup. While this was originally proven by \citep{Haettel} through different techniques, the approach via morphisms offers a particularly simple proof.

\medskip

Let \( r_1, \dots, r_{n+1} \) be the standard generators of the Artin group \(A[ B_{n+1}] \), where \( m_{r_n, r_{n+1}} = 4 \), and let \( t_0, \dots, t_n \) be the generators of the Euclidean braid group of type \(A[\widetilde{A}_n] \).
Define \( \rho := r_1 \cdots r_n r_{n+1} \). For \( 1 \leq i \leq n - 1 \), one checks that \( \rho r_i \rho^{-1} = r_{i+1} \). If we set \( r_0 := \rho r_n \rho^{-1} \), this relation extends cyclically modulo \( n+1 \), since \( \rho^2 r_n \rho^{-2} = r_1 \).

We now define an outer automorphism \( f \) of the group generated by \( t_0, \dots, t_n \) as follows:
\[
f: A[\widetilde{A}_n] \longrightarrow A[\widetilde{A}_n], \qquad t_i \longmapsto t_{i+1},
\]
where indices are taken modulo \( n+1 \). This gives an action of the infinite cyclic group \( \mathbb{Z} \cong \langle u \rangle \) on \( A[\widetilde{A}_n] \) by setting \( u \cdot g = u g u^{-1} := f(g) \). With this action, we form the semidirect product \( A[\widetilde{A}_n] \rtimes \langle u \rangle \), which has a presentation with generators \( \{t_0, \dots, t_n, u\} \) and relations consisting of the Artin relations for \( A[\widetilde{A}_n] \) together with
\[
ut_i u^{-1} = t_{i+1} \qquad \text{for } 0 \leq i \leq n \text{ (modulo } n+1).
\]

\begin{theorem}[\cite{KP02}]
The map
\[
\varphi: A[\widetilde{A}_n] \rtimes \langle u \rangle \longrightarrow A[B_{n+1}]
\]
defined by \( \varphi(t_i) = r_i \) and \( \varphi(u) = \rho \) is an isomorphism. In particular, the restriction of \( \varphi \) to \( A[\widetilde{A}_n] \) gives an embedding of \( A[\widetilde{A}_n] \) into \( A[B_{n+1}] \).
\end{theorem}

As a consequence of the definition of $\varphi$ we have the following lemma:

\begin{lemma}
Let \( \xi:A[B_{n+1}] \to \mathbb{Z} \) be the homomorphism defined by
\[
\xi(r_i) = 0 \quad \text{for } 1 \leq i \leq n, \qquad \xi(r_{n+1}) = 1.
\]
Then \( \xi(\rho) = 1 \), and the kernel of \( \xi \) is precisely \( \varphi(A[\widetilde{A}_n]) \).
\end{lemma}

We now prove the intersection result.

\begin{theorem}
Let \( P \) and \( Q \) be two parabolic subgroups of the Euclidean braid group \(A_S=A[ \widetilde{A}_n] \). Then \( P \cap Q \) is again a parabolic subgroup of \( A[\widetilde{A}_n ]\).
\end{theorem}

\begin{proof}
Since \( \varphi \) is an embedding, \( \varphi(P) \) (and \( \varphi(Q) \)) is a parabolic subgroup of \( A[B_{n+1}] \), even when \( P= g A_X g^{-1} \), $X\subset S$, and the set \( X \) includes \( t_0 \), since \( \varphi(t_0) = \rho r_n \rho^{-1} \) and  $P$ is proper. In that case, one can conjugate \( \varphi(P) \) inside \( \{r_1, \dots, r_n\} \) using a suitable power of \( \rho \).

\smallskip

Since the intersection of parabolic subgroups of \( A[B_{n+1}] \) is again a parabolic subgroup (\autoref{theorem:intersection}), we have that $\varphi(P \cap Q) = \varphi(P) \cap \varphi(Q)$ is a parabolic subgroup of \( A[B_{n+1}] \).

\smallskip

Now let \( P_2 = h A[B_{n+1}]_Y h^{-1} \) be a parabolic subgroup of \( A[B_{n+1}] \), with \( h \in A[B_{n+1}] \) and $Y\subset \{r_1,\dots, r_{n+1}\}$, and let \( P_1 = \varphi^{-1}(P_2) \subset A[\widetilde{A}_n] \). To finish the proof, we just need to show that \( P_1 \) is a parabolic subgroup of \( A[\widetilde{A}_n] \). First, observe that \( r_{n+1} \notin Y \): otherwise, \( hr_{n+1}h^{-1} \in P_2 \), hence in \( \varphi(A[\widetilde{A}_n]) \), which contradicts \( \xi(hr_{n+1}h^{-1}) = 1 \). Therefore, \( Y \subset \{r_1, \dots, r_n\} = \varphi(\{t_1, \dots, t_n\}) \).

\smallskip

Since \( \varphi \) is an isomorphism \( A[\widetilde{A}_n] \rtimes \langle u \rangle \to B_{n+1} \), we may write \( h = h_1 \rho^m \), with \( \varphi^{-1}(h_1) \in A[\widetilde{A}_n] \), \( m \in \mathbb{Z} \). Thus,
\[
P_2 = h_1 \rho^m A[B_{n+1}]_Y \rho^{-m} h_1^{-1},
\]
and using that \( u^m A_X u^{-m} = f^m(A_X) \), for $X\subset S$, we obtain:
\[
P_1 = \varphi^{-1}(P_2) = \varphi^{-1}(h_1) A_{f^m(\varphi^{-1}(Y))} \varphi^{-1}(h_1)^{-1},
\]
which is a parabolic subgroup of \( A[\widetilde{A}_n] \).
\end{proof}

\subsection{General results using retractions in the Salvetti complex}\label{section:retractionSalvetii}

\autoref{theorem_charneyparis} and \autoref{theorem_restandardisation} are two of the very few substantial results that hold for all Artin groups. Both are proven using a geometric retraction on the Salvetti complex. This section is dedicated to explaining this technique. Since we are interested in the use of the complex, we will skip some purely algebraic auxiliary  results that will be referenced. 

Let \( A_S = A_\Gamma \) be an Artin group with associated Coxeter group \( W_S \). Recall that \( S^f \) denotes the collection of subsets of \( S \) that define spherical-type parabolic subgroups. Also recall that the Salvetti complex \( \mathrm{Sal}(A_S) \) of \( A_S \) is the derived complex of the poset \( W_S \times S^f \), endowed with the partial order defined by: \( (u, X) \preceq (v, Y) \) if \( X \subset Y \), \( v^{-1}u \in W_Y \), and \( v^{-1}u \) is minimal in the coset \( v^{-1}uW_X \). The group \( W_S \) acts on the complex by sending \( (u, X) \mapsto (wu, X) \) for all \( w \in W_S \), so this action is by automorphisms. We also define \( \overline{\mathrm{Sal}}(A) \) to be the quotient of \( \mathrm{Sal}(A_S) \) by the action of \( W_S \).

\medskip

\begin{figure}[h!]
  \centering
  \def\svgwidth{1.2\columnwidth} 
  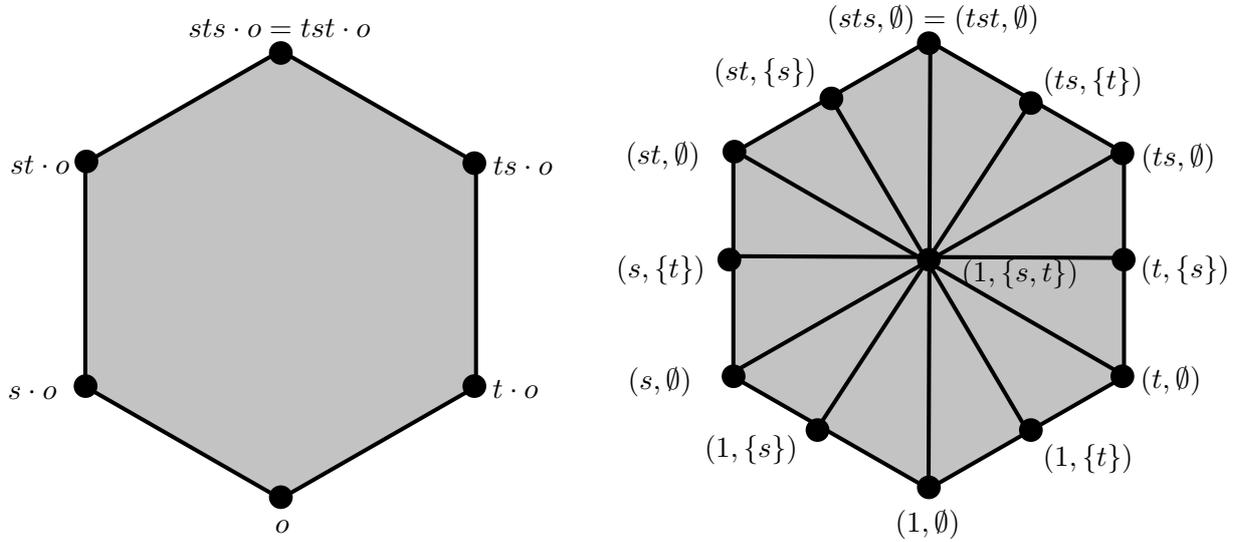
  \caption{On the left, a Coxeter cell corresponding to $W_\{{s,t\}}$ when $m_{s,t}=3$. On the right, $\mathbb{B}(1,\{s,t\})$ when $m_{s,t}=3$.}\label{Dibujo:CoxeterCell}
\end{figure}

Our first aim is to describe the cellular decompositions of both \( \mathrm{Sal}(A_S) \) and \( \overline{\mathrm{Sal}}(A) \). The building blocks of the Salvetti complex arise from the spherical-type parabolic subgroups. It is well known that any finite Coxeter group—such as \( W_X \) for \( X \in S^f \)—can be realized as a reflection group acting on \( \mathbb{R}^k \), where \( k = |X| \). Besides the classical references \citep{Bourbaki,Humphreys}, we refer the reader to Federica Gavazzi's PhD thesis \citep[Section~1.2]{GavazziThese} for a detailed exposition. We define the \emph{Coxeter cell} of \( W_X \) as the convex hull in \( \mathbb{R}^k \) of the orbit of a generic point \( o \), which is not fixed by any non-trivial element of \( W_X \).  For instance, if \( X = \{s, t\} \), then \( W_X \) is finite whenever \( m_{s,t} = m < \infty \), and in that case, the corresponding Coxeter cell is a regular \( 2m \)-gon (see \autoref{Dibujo:CoxeterCell}). Notice that this construction only makes sense when \( W_X \) is finite; otherwise, the orbit is infinite.

Next, we identify the Coxeter cells within the Salvetti complex. For any \( (u, X) \in W_S \times S^f \), we denote by \( \mathbb{B}(u, X) \) the subcomplex spanned by the set
\[
C(u, X) = \{(v, Y) \mid (v, Y) \preceq (u, X)\}.
\]
It is shown in \cite{ParisK} that \( \mathbb{B}(u, X) \) is homeomorphic to the Coxeter cell corresponding to~\( W_X \) (see again \autoref{Dibujo:CoxeterCell}). Thus, the Salvetti complex admits a cellular decomposition with cells \( \mathbb{B}(u, X) \) for all \( (u, X) \in W_S \times S^f \).

\begin{figure}[h!]
  \centering
  \def\svgwidth{0.5\columnwidth} 
  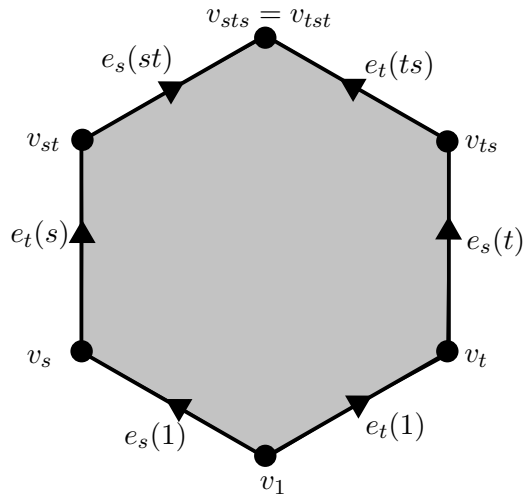
  \caption{Another notation for 2-skeleton of the Salvetti complex using our the cell decomposition.}\label{Dibujo:CoxeterCell2}
\end{figure}

In this decomposition we can see the vertices $v_w$ as $(w,\emptyset)$, where $w$ is word on two generators. Also, the edges $e_s(w)$ correspond to elements $(w,\{s\})$, where $s\in S$ and $w$ is a word in $s$ and other generator. If we orient the edges from $v_w$ to $v_{wt}$ where $l(wt)>l(w)$, we obtain \autoref{Dibujo:CoxeterCell2}. Notice that we can also have oriented edges in the opposite direction when $l(wt)<l(w)$ (this happens when $w$ can be written having $t$ as last letter). For example there is an edge $e_t(t)$ from~$v_t$ to $v_1$.

To see how this induces a cellular decomposition of $\overline{\mathrm{Sal}}(A)$, notice that the action of a non trivial element $w\in W_S$ sends $\mathbb{B}(u,X)$ to $\mathbb{B}(wu,X)$ in such a way that the the interior of the subcomplexes do not intersect, that is,
$$\mathrm{int}(\mathbb{B}(u,X))\cap\mathrm{int}(\mathbb{B}(wu,X))=\emptyset.$$
Then, when we quotient by $W$ we have that $\overline{\mathrm{Sal}}(A)$ have a cellular decomposition where each cell correspond to the orbit of $\mathbb{B}(1,X)$ for each $X\in S^f$.

\begin{exercise}\label{exo:Salvetti}
Prove that the 2-skeleton of $\overline{\mathrm{Sal}}(A)$ is the Cayley 2-complex of the standard presentation of $A_S$. That is, it has one vertex $x_0$, for each generator $s$ a directed loop $\overline{e}_s$, and for each relation $w=1$, a cell delimited by the edges that reads the word $w$. \emph{(See solution in page~\pageref{sol:Salvetti}.)}
\end{exercise}

Recall that $CA$ is the kernel of the natural epimorphism $\theta$ from the Artin group $A$ to the Coxeter group $W$. The short exact sequence  
$$1\longrightarrow CA \longrightarrow A \longrightarrow W \longrightarrow 1$$
corresponds to the regular covering \( \mathrm{Sal}(A_S) \rightarrow \overline{\mathrm{Sal}}(A) \), because $A=\pi_1(\overline{\mathrm{Sal}}(A),x_0)$ and  $CA=\pi_1({\mathrm{Sal}}(A),v_1)$ (by construction the loops in $\mathrm{Sal}(A)$ correspond to words representing the trivial element in $W_S$).

\medskip

Let $A_\Sigma$ be an Artin group with set of standard generators $\Sigma$ and let $W_S$ be its associated Coxeter group with set of generators $S$. We denote the generators differently to be able to differentiate words representing elements in the Artin or the Coxeter groups. The canonical surjection sends a standard generator $\sigma_i$ to a standard generator $s_i$.
Let $X\subset S$ and let $\Sigma_X\subset \Sigma$ the corresponding subset of generators of $\Sigma$. The main ingredient for our proofs will be a retraction described in \citep{GodelleParis}. First observe that the embedding $W_X\times X^f \hookrightarrow W_S\times S^f$ induces an embedding $i_X:\mathrm{Sal}(A_{\Sigma_X})\hookrightarrow \mathrm{Sal}(A_{\Sigma})$.

\begin{theorem}[{\citealp[Theorem~2.2]{GodelleParis}}]\label{theorem:retraction}
The embedding $i_X$ admits a retraction $\pi_X:\mathrm{Sal}(A_{\Sigma})\rightarrow \mathrm{Sal}(A_{\Sigma_X})$. This retraction is induced by the retraction 
$$\begin{array}{rrcl}
\pi'_X:&W_S\times S^f&\longrightarrow &W_X\times X^f \\
       &  (u,Y) & \longmapsto & (u_0,Y_0),
\end{array}
$$
where we have $u=u_0u_1$ with $u_0\in W_X$ and $u_1$ has minimal length in the coset $W_Xu_1$, and $Y_0 = X \cap u_1Yu_1^{-1}$.
\end{theorem}

In the above theorem, notice that $W_{Y_0}$ is finite because is a subset of $u_1W_Yu_1^{-1}$, which is a finite parabolic subgroup of $W_S$.

\begin{exercise}\label{exo:propiedades_retracto}
Let $W_S$ be a Coxeter group and $X\subset S$. Let $u=u_0u_1\in W_S$ with $u_0\in W_X$ and the length of $u_1$ is minimal in the coset $W_Xu_1$. Show that
\begin{enumerate}
\item $\pi_X(v_u)=v_{u_0}$.
\item For $s\in S$, let $x=u_1su_1^{-1}$. Then $\pi_X(e_s(u))=e_x(u_0)$ when $x\in X$ and otherwise $\pi_X(e_s(u))=v_{u_0}$. 
\end{enumerate}
\emph{(See solution in page~\pageref{sol:propiedadesretracto}.)}
\end{exercise}

Given a set of symbols $\mathcal{A}$, we denote by $\mathcal{A}^*$ the set of words that be written with those symbols. 

\begin{proposition}\label{prop:set_retraction}
Let $A_S$ be an Artin group with set of generators $\Sigma$. Let $W_S$ be its corresponding Coxeter group. Let $X\subset S$ and take the corresponding subset $\Sigma_X\subset \Sigma$. The retraction $\pi_X:\mathrm{Sal}(A_{S})\rightarrow \mathrm{Sal}(A_{X})$ induces a set retraction $$\widehat\pi_X:(\Sigma \sqcup \Sigma^{-1})^*\rightarrow (\Sigma_X \sqcup \Sigma_X^{-1})^*.$$
\end{proposition}

\begin{figure}[h]
  \centering
  \def\svgwidth{\columnwidth} 
  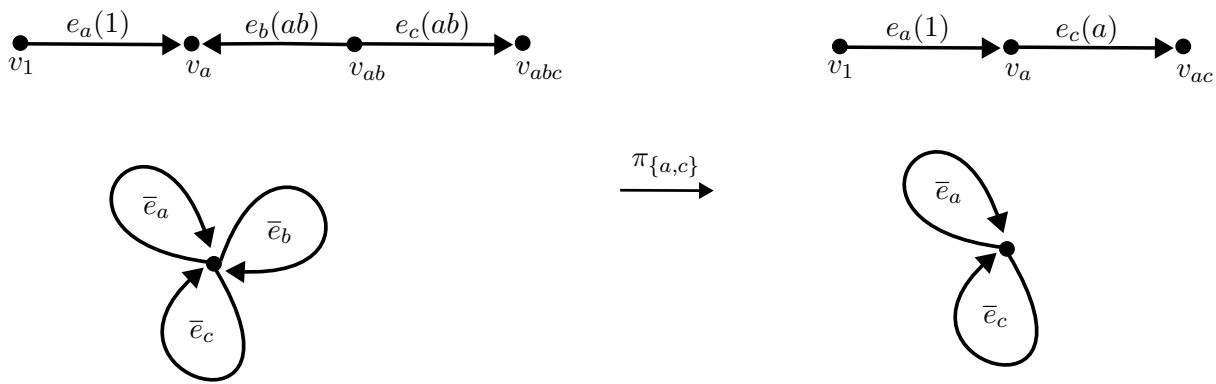
  \caption{We consider the Artin group $A_{\{a,b,c\}}$ with $m_{b,c}=2$ and we use the retract $\pi_\{a,c\}$. The word $ab^{-1}c$ is represented by a path with edges $e_a(1),e_b(ab),e_c(ab)$ in $\mathrm{Sal}(A_{\{a,b,c\}})$ whose images under the retract are, respectively, $e_a(1),v_a,e_c(a)$. In $\overline{\mathrm{Sal}}(A_{\{a,b,c\}})$ the first path corresponds to $\overline{e}_a,\overline{e}_b^{-1},\overline{e}_c$ and the second one to $\overline{e}_a,\overline{e}_c$, corresponding to the word $ac$ in $A_{\{a,b,c\}}$. }\label{Dibujo:ejemplo_Salvetti}
\end{figure}

\begin{proof}
We begin by considering a word 
\[
w = \sigma_{i_1}^{\epsilon_1} \cdots \sigma_{i_k}^{\epsilon_k}\in (\Sigma \sqcup \Sigma^{-1})^*,
\]
where each $\epsilon_j \in \{-1,1\}$. Using the retract we want to determine uniquely a word $w'\in (\Sigma_X \sqcup \Sigma_X^{-1})^*$. We will illustrate the proof with an example in \autoref{Dibujo:ejemplo_Salvetti}. By \autoref{exo:Salvetti}, we know that $A_\Sigma = \pi_1(\overline{\mathrm{Sal}}(A_\Sigma), x_0)$, so every generator in $\Sigma$ corresponds to a loop in $\overline{\mathrm{Sal}}(A)$. We then associate to the word $w$ the loop
\[
\overline{p}(w) = \overline{e}_{i_1}^{\epsilon_1} \cdots \overline{e}_{i_k}^{\epsilon_k}
\]
in $\overline{\mathrm{Sal}}(A)$.

Our first goal is to lift this loop to a path $p(w)$ in $\mathrm{Sal}(A)$ starting at the vertex $v_1$. To define this lift, we proceed as follows. For $\overline{e}_{i_1}^{\epsilon_1}$, we choose an edge $e_1$ starting at $v_1$. There are two such edges connecting $v_1$ to $v_{s_{i_1}}$: namely, $e_{s_{i_1}}(1)$, which we select if $\epsilon_1 = 1$, and $e_{s_{i_1}}(s_{i_1})$, which we select if $\epsilon_1 = -1$. Following this procedure, for $1 \leq j \leq k$, we define
\[
u_j = s_{i_1}^{\epsilon_1} \cdots s_{i_j}^{\epsilon_j} \in W,
\]
and set
\[
e_j = 
\begin{cases}
e_{s_{i_j}}(u_{j-1}) & \text{if } \epsilon_j = 1, \\
e_{s_{i_j}}(u_j) & \text{if } \epsilon_j = -1.
\end{cases}
\]
Thus, we define the lifted path as
\[
p(w) = e_1^{\epsilon_1} \cdots e_k^{\epsilon_k}.
\]

For the set retraction, we are only concerned with the 2-skeleton, so we use \autoref{exo:propiedades_retracto}. We decompose each $u_j$ as $u_j = u'_j u''_j$, where $u'_j \in W_X$ and $u''_j$ is the minimal-length representative of the coset $W_X u_j$. Define
\[
x_j = 
\begin{cases}
{u''}_{j-1} s_{i_j} {u''}_{j-1}^{-1} & \text{if } \epsilon_j = 1, \\[1mm]
{u''}_j s_{i_j} {u''}_j^{-1} & \text{if } \epsilon_j = -1.
\end{cases}
\]
Let $\sigma_{x_j}$ denote the Artin generator corresponding to $x_j\in S$, and set
\[
\chi_j = 
\begin{cases}
\sigma_{x_j} & \text{if } x_j \in X, \\
1 & \text{otherwise}.
\end{cases}
\]
We define
\[
\widehat\pi_X(w)=w' = \chi_1 \cdots \chi_k.
\]

It follows from \autoref{exo:propiedades_retracto} that $\pi_X(p(w))=p(w')$. Also, we have that $\pi_X(e_j)$ is the edge~$\overline{e}_{x_j}$ if $x_j \in X$, and the basepoint $x_0$ otherwise. Therefore, the image of $\pi_X(p(w))$ in $\overline{\mathrm{Sal}}(A)$ corresponds to the loop defined by the word $w'$. 
\end{proof}

\begin{proof}[Proof of \autoref{theorem_charneyparis}]
We want to prove that standard parabolic subgroups are convex with respect to the word length. Let $w=\sigma_{i_1}^{\epsilon_1} \cdots \sigma_{i_k}^{\epsilon_k}$ be a word representing an element of the Artin group $A_S$ generated by $\Sigma$. Let $W_S$ be the corresponding Coxeter group. Compute $w':=\widehat{\pi}_X(w)$ with the same notation as in the proof of \autoref{prop:set_retraction}. 
We prove two claims: 
\begin{enumerate}
\item If $w$ represent an element in $A_X$, then $w$ is equivalent to $w'$.
\item If $l(w)=l(w')$, then all the letters of $w$ belong to $\Sigma_X\sqcup \Sigma_X^{-1}$.
\end{enumerate}
Suppose that $w$ is a geodesic representing an element of $A_\Sigma$. By the first claim, $w'$ is equivalent to $w$. By construction, $l(w')\leq l(w)$ so $w$ is also geodesic and $l(w')\leq l(w)$, so by the second claim all the letters of $w$ belong to $\Sigma_X\sqcup \Sigma_X^{-1}$. This proves that $A_X$ is convex with respect to the word length. 

Let us now prove the claims. For the first one, choose a word $w''$ with letters in $\Sigma_X \sqcup \Sigma_X^{-1}$ that represents the same element as $w$. The corresponding loop $\overline{p}(w'')$ is homotopic to $\overline{p}(w)$ in $\overline{\mathrm{Sal}}(A)$ relative to the basepoint $x_0$. Let $p(w'')$ be the lift of $\overline{p}(v)$ starting at $v_0$ (as we did earlier in the proof for $p(w)$). Since $w''$ and $w$ represent the same group element, the paths $p(w'')$ and $p(w)$ have the same endpoints and are homotopic relative to those endpoints. Applying the retraction, we obtain that $\pi_X(p(w'')) = p(w'')$ and $\pi_X(p(w)) = p(w')$ are homotopic in $\mathrm{Sal}(A_X)$. Therefore, $w'$ also represents the same element as $w''$ (and hence as $w$).

For the second claim, suppose that $l(w) = l(w')$. Then, when applying the retraction to $w$, we must have $x_j \in X$ for every $j = 1, \dots, k$. We will prove by induction on $j$ that $s_{i_j} \in X$ for all $j = 1, \dots, k$.

For the base case $j=1$, we consider the two possibilities:
\begin{itemize}
    \item If $\epsilon_1 = 1$, then $x_1 = s_{i_1} \in X$, as desired.
    \item If $\epsilon_1 = -1$, suppose for contradiction that $s_{i_1} \notin X$. Then $s_{i_1}$ has minimal length in the coset $A_X s_{i_1}$, so $u_1 = u_1'' = s_{i_1}$ and hence $x_1 = s_{i_1} \in X$, which contradicts our assumption.
\end{itemize}
Now assume the claim holds for $1, \dots, j-1$, so that $u_{j-1} \in W_X$ and $u_{j-1}'' = 1$. We again distinguish two cases:
\begin{itemize}
    \item If $\epsilon_j = 1$, then $x_j = s_{i_j} \in X$.
    \item If $\epsilon_j = -1$, suppose $s_{i_j} \notin X$. Then, in the decomposition $u_j = u_j' u_j''$, we have $u_j' = u_{j-1}$ and $u_j'' = s_{i_j}$, so $x_j = s_{i_j} \in X$, again a contradiction.
\end{itemize}
This completes the inductive proof.
\end{proof}

\begin{exercise}\label{exo:palabrasSalvetti}
If two words $w$ and $w'$ represent the same element of $A_S$, then for every $X\subset S$ we have that $\widehat{\pi}_X(w)$ and $\widehat{\pi}_X(w')$ represent the same element in $A_X$. In other words, $\widehat\pi_X$ induces a set retraction $\widetilde{\pi}_X:A\rightarrow A_X$. \emph{(See solution in page~\pageref{sol:palabrasSalvetti}.)}
\end{exercise}

\begin{proof}[Proof of \autoref{theorem_restandardisation}]
Let $A = A_S$ be an Artin group, and let $P = \beta A_Z \beta^{-1}$ and $A_X$ be parabolic subgroups of $A_S$. We want to prove that there are $\alpha\in A_X$ and $Y\subset X$ such that $P=\alpha A_Y \alpha^{-1}$. 

As before, to differentiate sets of generators, we denote the standard set of generators of $A$ by $\Sigma=\{\sigma_1,\sigma_2,\dots\}$ in bijection with $S$. We aim to prove that there exist a subset $Y \subset X$ and an element $\alpha \in A_X$ such that
\[ P = \alpha A_Y \alpha^{-1}. \]
Let $g = \theta(\beta) \in W_S$. Then $gW_Zg^{-1} \subset W_X$. By standard Coxeter theory (see Lemma~2.2 in \cite{BlufsteinParis}), there exist $Y \subset X$ and $\beta_2 \in A_X$ such that
\[ \iota(g) A_Z \iota(g)^{-1} = \beta_2 A_Y \beta_2^{-1}. \]
Define $\beta_1 := \beta \iota(g)^{-1}$. Then we have:
\[
\beta A_Z \beta^{-1}
= \beta \iota(g)^{-1} \cdot \iota(g) A_Z \iota(g)^{-1} \cdot \iota(g) \beta^{-1}
= \beta_1 \beta_2 A_Y \beta_2^{-1} \beta_1^{-1}.
\]
Observe that both $\beta_2 A_Y \beta_2^{-1}$ and $\beta_1 \beta_2 A_Y \beta_2^{-1} \beta_1^{-1}$ lie in $A_X$, and that $\theta(\beta_1) = g g^{-1} = 1$, so $\beta_1 \in CA$. We now claim that for every $\gamma \in A_X$, we have
\[ \beta_1 \gamma \beta_1^{-1} = \pi_X(\beta_1) \gamma \pi_X(\beta_1)^{-1}. \]
If this claim holds, then
\[
\beta A_Z \beta^{-1} = \beta_1 \beta_2 A_Y \beta_2^{-1} \beta_1^{-1}
= \pi_X(\beta_1) \beta_2 A_Y \beta_2^{-1} \pi_X(\beta_1)^{-1},
\]
so we may take $\alpha := \pi_X(\beta_1) \beta_2 \in A_X$, completing the proof.

Let us now prove the claim. Consider a word
\[
w_1=\sigma_{i_1}^{\epsilon_1}\sigma_{i_2}^{\epsilon_2}\cdots \sigma_{i_k}^{\epsilon_k}
\]
representing $\gamma$, where all the letters belong to $\Sigma_X \sqcup \Sigma_X^{-1}$, and a word
\[
w_2=\sigma_{j_1}^{\mu_1}\sigma_{j_2}^{\mu_2}\cdots \sigma_{j_{k'}}^{\mu_{k'}}
\]
representing $\beta_1$, with letters in $\Sigma \sqcup \Sigma^{-1}$. To compute $\widetilde\pi_X(\beta_1 \gamma \beta_1^{-1})$, we use the proof of \autoref{prop:set_retraction}. We consider the word
\[
w = \sigma_{j_1}^{\mu_1}\sigma_{j_2}^{\mu_2}\cdots \sigma_{j_{k'}}^{\mu_{k'}} 
\left(\sigma_{i_1}^{\epsilon_1}\sigma_{i_2}^{\epsilon_2}\cdots \sigma_{i_k}^{\epsilon_k}\right)
\sigma_{j_{k'}}^{-\mu_{k'}}\sigma_{j_{k'-1}}^{-\mu_{k'-1}}\cdots \sigma_{j_1}^{-\mu_1}.
\]
We define:
\[
\begin{aligned}
u_{q,1} &=
\begin{cases}
1 & \text{if } q = 0, \\
s_{j_1} s_{j_2} \cdots s_{j_q} & \text{if } 1 \leq q \leq k',
\end{cases}
\\[2mm]
u_{q,2} &=
\begin{cases}
\theta(\beta_1) = 1 & \text{if } q = 0, \\
s_{i_1} s_{i_2} \cdots s_{i_q} & \text{if } 1 \leq q \leq k,
\end{cases}
\\[2mm]
u_{q,3} &=
\begin{cases}
\theta(\beta_1)\theta(\alpha) = \theta(\alpha) & \text{if } q = 0, \\
\theta(\alpha)\, s_{j_{k'}} s_{j_{k'-1}} \cdots s_{j_{k'-q+1}} & \text{if } 1 \leq q \leq k'.
\end{cases}
\end{aligned}
\]

We now explain these definitions in more detail. The sequence \( u_{q,1} \) corresponds to the initial segment
$
\sigma_{j_1}^{\mu_1} \sigma_{j_2}^{\mu_2} \cdots \sigma_{j_q}^{\mu_q},
$
the sequence \( u_{q,2} \) corresponds to the central part
$
\sigma_{i_1}^{\epsilon_1} \sigma_{i_2}^{\epsilon_2} \cdots \sigma_{i_q}^{\epsilon_q},
$
and \( u_{q,3} \) corresponds to the final segment
$
\sigma_{j_{k'}}^{-\mu_{k'}} \sigma_{j_{k'-1}}^{-\mu_{k'-1}} \cdots \sigma_{j_{k'-q+1}}^{-\mu_{k'-q+1}}.
$
When defining \( u_{0,2} \), we must take into account the position in the Salvetti complex at which the middle segment begins. In the Salvetti complex, vertices correspond to elements of the Coxeter group. If we choose a path representing \( w \) starting at \( v_0 \), then after reading the prefix \( \sigma_{j_1}^{\mu_1} \cdots \sigma_{j_{k'}}^{\mu_{k'}} \), we reach the vertex corresponding to \( \theta(\beta_1) = 1 \). Similarly, after reading the central segment \( \sigma_{i_1}^{\epsilon_1} \cdots \sigma_{i_k}^{\epsilon_k} \), we reach the vertex \( \theta(\beta_1)\theta(\alpha) = \theta(\alpha) \).

Next, we write \( u_{q,r} = u_{q,r}' u_{q,r}'' \), for \( r = 1,2,3 \), where \( u_{q,r}' \in A_X \) and the length of \( u_{q,r}'' \) is minimal in the coset \( W_X u_{q,r}'' \), and we define:

\[
\begin{array}{ll}
\begin{array}{rcl}
x_{q,1} &=&
\begin{cases}
u_{q-1,1}''\, s_{j_q}\, u_{q-1,1}''^{-1} & \text{if } \epsilon_q = 1 \\
u_{q,1}''\, s_{j_q}\, u_{q,1}''^{-1}     & \text{if } \epsilon_q = -1
\end{cases}
\end{array},
&
1\leq q\leq k'
\\[2 mm]
\begin{array}{rcl}
x_{q,2} &=&
\begin{cases}
u_{q-1,2}''\, s_{i_q}\, {u_{q-1,2}''}^{-1} & \text{if } \epsilon_q = 1 \\
u_{q,2}''\, s_{j_q}\, {u_{q,2}''}^{-1}     & \text{if } \epsilon_q = -1
\end{cases}
\end{array},
&
1\leq q\leq k
\\[2 mm]
\begin{array}{rcl}
x_{q,3} &=&
\begin{cases}
u_{q-1,3}''\, s_{i_{k'-q+1}}\, u_{q-1,3}''^{-1} & \text{if } \epsilon_{k'-q+1} = 1, \\
u_{q,3}''\, s_{j_{k'-q+1}}\, u_{q,3}''^{-1}     & \text{if } \epsilon_{k'-q+1} = -1
\end{cases}
\end{array},
&
1\leq q\leq k'
\end{array}
\]
And finally we denote by $\sigma_{x_{q,r}}$ the standard Artin generator that correspond to $x_{q,r}$ and we define
$$
\begin{array}{rcl}
\chi_{q,r} &=&
\begin{cases}
\sigma_{x_{q,r}} & \text{if } x_{q,r}\in X \\
1 & \text{otherwise.} 
\end{cases}
\end{array}
$$
Then the retracted word is 
\[
\widehat\pi_X(w) = w' = \chi_{1,1} \cdots \chi_{k',1}\chi_{1,2} \cdots \chi_{k,2}\chi_{1,3}\cdots \chi_{k',3} = \widehat\pi_X(w_2)\widehat\pi_X(w_1)\chi_{1,3}\cdots \chi_{k',3}.
\]
Notice that, since all the $\sigma_{j_q}$ lie in $\Sigma_X$, we have that $\pi_X(w_1)=w_1$. We also claim that $\widehat\pi_X(w_2)^{-1}=\chi_{1,3}\cdots \chi_{k',3}$. If this is true, since $\beta_1\alpha\beta_1^{-1}\in A_X$, by \autoref{exo:palabrasSalvetti} we will have
$$\beta_1\alpha\beta_1^{-1}=\widetilde\pi_X(\beta_1\alpha\beta_1^{-1})=\widetilde\pi_X(\beta_1)\alpha\widetilde\pi_X(\beta_1)^{-1}.$$
To prove that $\widehat\pi_X(w_2)^{-1}=\chi_{1,3}\cdots \chi_{k',3}$, we need to prove that $\chi_{q,1}=\chi_{k'-q+1,3}^{-1}$ for $1\leq q \leq k'$. Notice that 
$1=\theta(\beta_1)=s_{j_1} s_{j_2} \cdots s_{j_k'}$, so we have
$$s_{j_1} s_{j_2} \cdots s_{j_q}= s_{j_{k'}} s_{j_{k'-1}} \cdots s_{j_{q+1}}.$$
Notice that, since $\theta(\alpha)\in W_X$ and $s_{j_q}^{-1}=s_{j_q}$, we have that $u''_{q,3}=u''_{k'-q,1}$. Then, for $1\leq q\leq k'$, we have
$$\begin{array}{rcl}
x_{q,3} &=&
\begin{cases}
u_{k'-q+1,3}''\, s_{i_{k'-q+1}}\, u_{k'-q+1,3}''^{-1}& \text{if } \epsilon_{k'-q+1} = 1, \\
u_{k'-q,3}''\, s_{j_{k'-q+1}}\, u_{k'-q,3}''^{-1} & \text{if } \epsilon_{k'-q+1} = -1
\end{cases}
\end{array}.$$
In both cases we obtain $x_{q,3}=x_{k'-q+1}^{-1}$, and this implies $\chi_{q,1}=\chi_{k'-q+1,3}^{-1}$ as we wanted.
\end{proof}

\section{Solutions to the exercises}

\noindent
\textbf{\autoref{exo:GarsideElement}}. \emph{Compute the Garside element of $A[H_3]$.} 

\begin{solution}\label{sol:Garside}
We now construct a positive word \( w \) representing \( \Delta \).  
Recall that
\[
\Delta = \bigvee \{a,b,c\},
\]
where
\[
m_{a,c} = 2, \quad m_{a,b} = 3, \quad m_{b,c} = 5.
\]

We will parenthesize a subword wherever we are going to apply a relation. The least common multiple of \(a\) and \(c\) is simply \(ac = ca\).  
We next append a positive suffix so that \(b\) also appears as a prefix.  
To insert \(b\), we add \(bcbc\), obtaining:
\[
a(cbcbc) = abcbcb.
\]
However, this is still insufficient, since to have \(b\) as a prefix we must obtain \(aba\) as a prefix.  
Appending \(ab\) yields:
\[
abcbc(bab) = abcb(ca)ba = abc(ba)cba.
\]
The parenthesized term must be completed with a \(b\).  
Thus, we add \(c\) to \(cba\):
\[
abcbacb(ac) = abcba(cbc)a.
\]
The parenthesized term now needs \(bc\) to be completed.  
To obtain a \(b\), we add \(ba\):
\[
abcbacbc(aba) = abcba(cbcb)ab.
\]
We still require a \(c\) to complete the parentheses, so we add \(cbcb\):
\[
abcbacbcba(bcbcb) 
= abcbacbcb(ac)bcbc 
= abcba(cbcbc)abcbc 
= abc(bab)cbcbabcbc
\]
\[
= ab(ca)bacbcbabcbc 
= (aba)cbacbcbabcbc 
= babcbacbcbabcbc.
\]
Therefore,
\[
\Delta = babcbacbcbabcbc.
\]
Observe that this word has \(15\) letters.  
Moreover, it is equivalent to \((abc)^5\), since
\[
(bab)cbacbc(bab)cbc = ab(ac)b(ac)bcab(ac)bc = abcabcabcabcabc.
\]

\end{solution}

\noindent
\textbf{\autoref{exo_inf_sum_mixednoralforms}}. \emph{Let $\alpha=a^{-1}b$ be in mixed np-normal form. Prove that, if $p<0$, $\mathrm{sup}(a)=-\mathrm{inf}(\alpha)$ and $\mathrm{sup}(b)=\mathrm{sup}(\alpha)$.}
\begin{solution}\label{sol:inf}Consider the left normal form of \(\alpha\):
\[
\alpha = \Delta^p x_1 \cdots x_r.
\]
We define
\[
a = \big( \Delta^p x_1 \cdots x_{-p} \big)^{-1},
\quad
b = x_{-p+1} \cdots x_r.
\]

First, observe that
\[
\Delta^p x_1 \cdots x_{-p}
= \Delta^{-1} \big( \Delta^{p+1} x_1 \Delta^{-(p+1)} \big) 
\cdot \Delta^{-1} \big( \Delta^{p+2} x_2 \Delta^{-(p+2)} \big) 
\cdots \Delta^{-1} x_{-p}.
\]
By the properties of the Garside structure, each term
\(\Delta^{p+i} x_i \Delta^{-(p+i)}\) is a simple element;  
multiplication by \(\Delta^{-1}\) therefore yields a negative element.  
Hence \(a\) is totally positive.

It remains to show that \(a \wedge b = 1\).  
By \autoref{Prop:negativo}, the left normal form of \(a\) is
\[
a =  
\big( \Delta^{-1} x'_{-p} \Delta \big)
\cdot \big( \Delta^{-2} x'_{-p-1} \Delta^{2} \big)
\cdots
\big( \Delta^{p} x'_{1} \Delta^{-p} \big),
\]
where \(x'_i = \Delta x_i^{-1}\).  
The element \(b\) is already in left normal form, so it has no \(\Delta\) as a prefix.  
Thus, the only possible common prefixes letters of \(a\) and \(b\) are letters that are simultaneously prefixes of  
\(\Delta^{-1} x'_{-p} \Delta = x_{-p}^{-1} \Delta\) and \(x_{-p+1}\).  

Suppose \(s\) is such a letter.  
Since \(\Delta = x_{-p} \cdot x_{-p}^{-1} \Delta\), if \(s\) were a prefix of \(x_{-p}^{-1} \Delta\) then \(x_{-p} s\) would be a simple element.  
This contradicts the fact that \(x_{-p} \cdot x_{-p+1}\) is in left normal form.  
Therefore, \(a \wedge b = 1\).

\end{solution}

\noindent
\textbf{\autoref{exo_parabolicclosure}}. \emph{Prove that for every parabolic subgroup $Q$, we have that $P_{z_Q}=Q$.}
\begin{solution}\label{sol:parabolicclosure}
Let \(A_S\) be an Artin group of spherical type. We can write 
\[
Q = \alpha A_X \alpha^{-1}, \quad z_Q = \alpha \Delta_X \alpha^{-1},
\]
for some \(X \subset S\) and \(\alpha \in A_S\).  
As shown in the proof of \autoref{theorem:unique_parabolic_closure}, we have  
\[
P_{\alpha \Delta_X \alpha^{-1}} = \alpha\, P_{\Delta_X}\, \alpha^{-1}.
\]
Since \(\Delta_X\) is positive, it is already a recurrent element. Also, the support of $\Delta_X$ is \(X\), so \(P_{\Delta_X} = A_X\).  
Therefore,
\[
P_{z_Q} = \alpha A_X \alpha^{-1}.
\]
\end{solution}

\noindent
\textbf{\autoref{exo_recurrent_swapping}}. \emph{Prove that when we apply recurrent swapping to an element, the infimum can only increase and the supremum can only decrease.}
\begin{solution}\label{sol:recurrent}
Consider the left normal form of \(\alpha\):
\[
\alpha = \Delta^p x_1 \cdots x_r.
\]
If \(p \geq 0\), then the mixed normal form of \(\alpha\) is \(1 \cdot \alpha\), and applying the swap operation leaves the element unchanged.  
If \(p < 0\), as seen in the solution of \autoref{exo_inf_sum_mixednoralforms}, the mixed normal form of \(\alpha\) is \(a^{-1} b\) with
\[
a = \big( \Delta^p x_1 \cdots x_{-p} \big)^{-1},
\quad
b = x_{-p+1} \cdots x_r.
\]
After applying a swap, we obtain
\[
x_{-p+1} \cdots x_r \, \Delta^p x_1 \cdots x_{-p}
= \Delta^{p} \, \Delta^{-p} x_{-p+1} \cdots x_r \, \Delta^{p}  x_1 \cdots x_{-p}.
\]
Since the conjugate of a simple element by a power of \(\Delta\) is still a simple element, this expression consists of a power of \(\Delta\) followed by a product of simple elements \(y_1 \cdots y_r\).  

However, it may happen that \(y_i \cdot y_{i+1}\) is not left-weighted. In this case, we compute
\[
y_i' = y_i y_{i+1} \wedge \Delta,
\quad
y_{i+1}' = {y_i'}^{-1} y_i y_{i+1},
\]
so that \(y_i' \cdot y_{i+1}'\) is left-weighted.  
This operation is called a \emph{left sliding}. Applying it iteratively yields a (non-optimized) algorithm to obtain a normal form.  

During this process, it may happen that some \(y_i y_{i+1}\) has \(\Delta\) as a prefix, which increases~\(p\).  
Since all simple elements are positive, \(p\) can never decrease.  
It may also happen that \(y_i y_{i+1} \wedge \Delta = y_i y_{i+1}\), which reduces the number of simple factors in the normal form.  
Since left sliding always produces at most two non-trivial simple elements, the number of simple factors can never increase.
\end{solution}
%

\noindent
\textbf{\autoref{exercise:pointwise_geodesic}}. \emph{If \( g \in A_S\) fixes two vertices of the cubical decomposition of the Deligne complex of an FC-type Artin group $A_S$, then it fixes pointwise any combinatorial geodesic between them.}
\begin{solution}[\citealp[Remark 2.1]{Rose}]\label{sol:pointwise}
If a point $x$ in the interior of a cube is fixed by the action of an element $g$, then $g$ fixes the entire cube. This follows because each vertex of a cube corresponds to a coset of a parabolic subgroup $A_X$ for some $X \subset S$, and the action of $A_S$ sends such a coset to another coset of $A_X$. Since no two vertices of a cube are cosets of the same subgroup $A_X$, there is no non-trivial element that can permute the vertices of the cube.  

Now, suppose $g$ fixes two vertices $v_1$ and $v_2$. Then $g$ fixes not only every minimal length edge path between $v_1$ and $v_2$, but also the whole interval of cubes containing such a path. In particular, $g$ fixes pointwise any combinatorial geodesic between $v_1$ and $v_2$.
\end{solution}

\noindent
\textbf{\autoref{exercise:simplices}}. \emph{Let $A$ be an Artin group and $K$ be its Artin complex. Prove that if \( P \) and \( Q \) are parabolic subgroups, corresponding to the stabilizers of simplices \( D_1 \) and \( D_2 \) of $K$ that share a face \( F \), then, up to conjugacy, \( P \) and \( Q \) are parabolic subgroups of \( \mathrm{Stab}(F) \).}
\begin{solution}\label{sol:simplices}
By construction, since $F\subseteq D_1$ and $F\subseteq D_2$, we have that $\mathrm{Stab}(D_1)\subseteq\mathrm{Stab}(F)$ and  $\mathrm{Stab}(D_1)\subseteq\mathrm{Stab}(F)$. Now, up to conjugacy, we can suppose that~$\mathrm{Stab}(F)$ is standard, so an Artin group with standard generators. Using \autoref{theorem_restandardisation}, we have that~$P$ and~$Q$ are parabolic subgroups of~$\mathrm{Stab}(F)$.
\end{solution}

\noindent
\textbf{\autoref{exo:FC-decomposition}}. \emph{Let \( A_\Gamma = A_S \) be an Artin group with associated Coxeter graph \( \Gamma \) (with the no-$\infty$ convention), and let \( \Gamma_1 \) and~\( \Gamma_2 \) be two subgraphs of \( \Gamma \) such that \( \Gamma = \Gamma_1 \cup \Gamma_2 \) and \( \Gamma_1 \cap \Gamma_2 \neq \emptyset \). Show that \( A_\Gamma \) admits the following amalgamated product decomposition:
\[
A_\Gamma = A_{\Gamma_1} *_{A_{\Gamma_1 \cap \Gamma_2}} A_{\Gamma_2}.
\]
In particular, if \( m_{s,t} = \infty \) for some \( s,t \in S \), then:
\[
A_S = A_{S \setminus \{s\}} *_{A_{S \setminus \{s,t\}}} A_{S \setminus \{t\}}.
\]
This shows that if \( A_S \) is of FC type, then it can be expressed as an amalgamated product of spherical-type Artin groups.}
\begin{solution}\label{sol:FC}
For the first part, observe that the standard presentations of 
\(A_{\Gamma_1}\) and \(A_{\Gamma_2}\) share precisely the relations coming from 
the standard presentation of \(A_{\Gamma_1 \cap \Gamma_2}\). 
This matches the definition of an amalgamated product, with the injective homomorphisms 
given by the natural inclusions. 
For the second part, note that any subset \(X \subset S\) free of \(\infty\) 
generates a standard parabolic subgroup of spherical type.
\end{solution}

\noindent
\textbf{\autoref{exo:retracto}}. \emph{Consider the Artin group \( A_S \) corresponding to the graph in page \pageref{grafotikz}. Prove that if \( X = \{a, c, d, e, f\} \) and \( Y = \{b, c, d\} \), then for any \( g, h \in A_S \), the intersection \( gA_Xg^{-1} \cap hA_Yh^{-1} \) is a parabolic subgroup over \( \{c, d\} \), contained in \( A_{\{c, d, e, f\}} \).}
\begin{solution}\label{sol:retracto}
Suppose \( X = \{a, c, d, e, f\} \) and \( Y = \{b, c, d\} \). Up to conjugation by \( h \), we may assume \( h = 1 \), so that \( Q = A_Y \), and \( g \in A_{\{a, b\}} \). Applying the retraction \( \rho_Y \), we obtain:
\[
P \cap Q = P \cap A_Y = \rho_Y(P \cap A_Y) \subseteq \rho_Y(g)A_{X \cap Y}\rho_Y(g)^{-1} = A_{X \cap Y} \subseteq A_{X'}.
\]
Note that \( X' = \{c, d, e, f\} \) is a direct component of \( X \), and that \( P \cap A_{X'} = A_{X'} \), so:
\[
P \cap Q = P \cap A_{X'} \cap A_Y = A_{X'} \cap A_Y = A_{\{c, d\}}.
\]
Thus, \( P \cap Q \) is a parabolic subgroup over \( X \cap Y \), contained in \( X' \), as claimed.
\end{solution}

\noindent
\textbf{\autoref{exo:Salvetti}}. \emph{Prove that the 2-skeleton of $\overline{\mathrm{Sal}}(A)$ is the Cayley 2-complex of the standard presentation of $A_S$. That is, it has one vertex $x_0$, for each generator $s$ a directed loop $\overline{e}_s$, and for each relation $w=1$, a cell delimited by the edges that reads the word $w$.}
\begin{solution}\label{sol:Salvetti}
Since the vertices \( v_w \) are in bijective correspondence with the elements of \( W \), they all lie in the same equivalence class when taking the quotient by \( W \). Hence, \(\overline{\mathrm{Sal}}(A)\) has a single vertex \(x_0\). Similarly, any two edges \( e_s(w) \) and \( e_s(w') \) belong to the same equivalence class, so in~\(\overline{\mathrm{Sal}}(A)\) there is one loop \(\overline{e}_s\) for each generator \(s\). Finally, we have seen that the relations \( w=1 \) of the standard presentation of the Artin group correspond to traversing the boundary of the 2-cells, as illustrated in \autoref{Dibujo:CoxeterCell2}. Thus, in \(\overline{\mathrm{Sal}}(A)\) we obtain one 2-cell bounded by \(w'\) for each relation \(w'=1\) in the group.
\end{solution}

\noindent
\textbf{\autoref{exo:propiedades_retracto}}. \emph{Let $W_S$ be a Coxeter group and $X\subset S$. Let $u=u_0u_1\in W_S$ with $u_0\in W_X$ and the length of $u_1$ is minimal in the coset $W_Xu_1$. Show that
\begin{enumerate}
\item $\pi_X(v_u)=v_{u_0}$.
\item For $s\in S$, let $x=u_1su_1^{-1}$. Then $\pi_X(e_s(u))=e_x(u_0)$ when $x\in X$ and otherwise $\pi_X(e_s(u))=v_{u_0}$. 
\end{enumerate}}
\begin{solution}\label{sol:propiedadesretracto}
We already know that the retraction is well defined. Moreover, a retraction is always a continuous map.  
The vertex \( v_u \) corresponds to \((u,\emptyset)\), which is sent by the retraction to \((u_0,\emptyset)\). In other words, \(v_u\) is mapped to \(v_{u_0}\).  
Now, the edge \(e_s(u)\) corresponds to \((u,\{s\})\), which is sent by the retraction to \((u_0,Y_0)\), where 
\[
Y_0 = X \cap u_1 \{s\} u_1^{-1}.
\]  
If \(s \in X\), then \(Y_0 = \{s\}\), and consequently \(e_s(u)\) is mapped to \(e_s(u_0)\). Otherwise, if \(s \notin X\), then \(Y_0 = \emptyset\) and \(e_s(u)\) is mapped to the vertex \(v_{u_0}\).
\end{solution}

\noindent
\textbf{\autoref{exo:palabrasSalvetti}}. \emph{If two words $w$ and $w'$ represent the same element of $A_S$, then for every $X\subset S$ we have that $\widehat{\pi}_X(w)$ and $\widehat{\pi}_X(w')$ represent the same element in $A_X$. In other words, $\widehat\pi_X$ induces a set retraction $\widetilde{\pi}_X:A\rightarrow A_X$.}
\begin{solution}[{\citealp[Proposition~2.3(1)]{BlufsteinParis}}]\label{sol:palabrasSalvetti}
Since \( w \) and \( w' \) are equivalent, the corresponding loops 
\(\overline{p}(w)\) and \(\overline{p}(w')\) in 
\(\overline{\mathrm{Sal}}(A_S)\) represent the same element, and therefore 
they are homotopic. We take lifts \(p(w)\) and \(p(w')\) starting at \(v_1\). 
As they represent the same element, the lifts also end at the same vertex and 
are homotopic relative to their endpoints. Since a retraction is continuous, 
the images \(p_X(w):=\pi_X(p(w))\) and \(p_X(w'):=\pi_X(p(w'))\) are also 
homotopic relative to their endpoints. Projecting again to 
\(\overline{\mathrm{Sal}}(A_X)\), we obtain two loops~\(\overline{p}_X(w)\) 
and~\(\overline{p}_X(w')\) that represent the same element in~\(A_X\). 
Finally, as shown in the proof of \autoref{prop:set_retraction}, the words 
\(\widehat{\pi}_X(w)\) and \(\widehat{\pi}_X(w')\) correspond precisely to 
those loops, and hence they represent the same element of \(A_X\), as desired.
\end{solution}

\bigskip

\noindent{\textbf{\Large{Acknowledgments}}}

\smallskip

I thanks the organisers of Winterbraids XIV for giving me the opportunity of speak and write about this topic. I was supported bt the research project PID2022-138719NA-I00, financed by MCIN/AEI/10.13039/501100011033/FEDER, UE, and by a Ram\'on y Cajal 2021 grant, also financed by the Spanish Ministry of Science and Innovation. 
I am grateful to the anonymous referee for their careful reading and helpful suggestions. I would also like to thank my great-aunt, Eloísa Pastor Ramírez, who lovingly hosted me while I was writing a substantial part of this manuscript.

\medskip

\medskip
\bibliography{bibliography}

\bigskip\bigskip{\footnotesize%

\noindent
\textit{\textbf{María Cumplido} \\ 
Departmento de \'Algebra,
Facultad de Matem\'aticas,
Universidad de Sevilla, \\
and Instituto de Matemáticas de la Universidad de Sevilla (IMUS).\\
Calle Tarfia s/n
41012, Seville, Spain.} \par
 \textit{E-mail address:} \texttt{\href{mailto:cumplido@us.es}{cumplido@us.es}}
 }

\end{document}